# An adaptive phase field framework for large-scale interface evolution problems using a strong-form gradient smoothing approach

Zirui Mao[1*], Alice Xu[2], Shenyang Hu[1], Panos Stinis[1], Yuyan Shao[1], Ang Li[1], Yulan Li[1]

[1] *Pacific Northwest National Laboratory, Richland, WA 99352*
[2]*Hanford High School, Richland, WA 99352*

E-mail: zirui.mao@pnnl.gov;

**Abstract:**

Multiscale problems with evolving interfaces are ubiquitous in science and engineering. Phase-field models are a powerful tool for simulating interface-dominated phenomena in computational mechanics and materials modeling, but their application to large-scale problems is often constrained by the high computational cost of resolving thin diffuse interfaces over the entire domain. This paper presents an efficient strong-form phase-field solver that couples the Gradient Smoothing Method (GSM) with a hierarchical adaptive and moving structured mesh, enabling automatic localization of resolution within a narrow interfacial region while retaining coarse discretization in bulk domains. A layered refinement design is introduced to preserve locally uniform resolution across the interface, allowing the GSM discretization to maintain overall second-order accuracy despite strong mesh non-uniformity away from the interface. Although GSM incurs a higher per-degree-of-freedom cost than standard finite-difference schemes, the adaptive framework substantially reduces the total number of degrees of freedom, resulting in near-linear computational scaling compared with the quadratic scaling of uniform-grid approaches. Numerical examples based on the Allen–Cahn and Cahn–Hilliard equations demonstrate that the proposed adaptive GSM solver delivers excellent accuracy for interface evolution while attaining more favorable computational complexity, $O(N)$, than existing weak-form and strong-form solvers, becoming significantly more efficient for large-scale problems with thin interfaces or a small interfacial area fraction relative to the whole domain.

## 1. Introduction

Multiscale problems with evolving interfaces arise across fracture and crack propagation [1-4], fluid mechanics [5, 6], and microstructure evolution [7, 8], where sharp boundaries between distinct phases migrate, merge, and undergo topological change. The phase-field method has emerged as a standard modeling framework for such problems by replacing sharp interfaces with a diffuse transition layer described by a continuous order parameter [7, 9]. This formulation enables robust representation of complex interface kinetics and topology changes, and naturally accommodates coupling with additional physics, such as diffusion, elasticity, and fluid flow. Despite these advantages, phase-field methods face a persistent computational bottleneck: accurate resolution of the thin diffuse interface requires a spatial mesh fine enough to capture steep gradients across the interfacial band. Applying such resolutions uniformly over the entire domain leads to a rapid increase in degrees of freedom, memory footprint, and computational cost, particularly in large-scale or three-dimensional simulations. This inefficiency strongly motivates numerical strategies that dynamically concentrate resolution in regions of interfacial activity while maintaining a coarser discretization in the bulk phases, where solution fields vary smoothly and fine resolution is unnecessary.

Over the last decades, two complementary strands of numerical technology have matured that are directly relevant to efficient phase-field simulations with evolving interfaces. On one hand, adaptive mesh refinement (AMR) [10-12] strategies within weak-form finite-element frameworks [1] have been widely adopted for phase-field modeling. In these approaches, structured meshes are dynamically refined in the vicinity of diffuse interfaces based on indicators such as phase-field gradients, energy density, or distance measures to the interface, enabling substantial reductions in computational cost by concentrating degrees of freedom in localized regions of interface. Representative examples include adaptive FEM implementations in general multiphysics frameworks such as MOOSE [13], FEniCS [3, 14], Alamo [15], and PRISMS-PF [16], where structured AMR is employed to track evolving interfaces and cracks. While these weak-form AMR approaches offer high geometric flexibility and broad multiphysics capability, they typically rely on repeated global matrix assembly and the construction of large Jacobian systems, making the numerical solvers highly sensitive to mesh quality and to the degree of nonlinearity in the governing equations. For highly nonlinear phase-field problems, the resulting discrete systems can become very large and poorly conditioned, thereby making iterative solvers difficult to converge and, in extreme cases, prone to failure. These issues can ultimately limit scalability and efficiency on modern high-performance computing architectures.

On the other hand, strong-form, stencil-based discretizations (such as generalized finite-difference [17, 18] and finite-volume methods [2, 19, 20] with block-structured AMR) discretize the governing phase-field equations directly on structured or semi-structured node sets. In the context of interface evolution, strong-form solvers have successfully employed block-structured or narrow-band AMR to refine only a thin region surrounding the diffuse interface while keeping the bulk coarse. For instance, adaptive finite-difference solvers have been used for dendritic solidification and alloy solidification by refining grids in regions of steep phase-field gradients [21], while more recent block-structured AMR frameworks solve coupled phase-field and elasticity equations in a matrix-free manner using local stencils and ghost layers [22]. These approaches enable deep local refinement, automatic remeshing, and efficient parallel execution, while avoiding global stiffness assembly. Gradient-smoothing ideas naturally extend this strong-form philosophy: they have been explored both within smoothed finite-element formulations [23] and as standalone strong-form schemes [24-26] that generalize classical finite-difference operators to irregular or non-uniform node distributions, thereby improving the robustness and accuracy of derivatives approximation on distorted or adaptively refined meshes.

Compared with weak-form approaches, strong-form solvers offer several practical advantages that are particularly attractive for nonlinear phase-field dynamics. Since the differential operators are discretized directly, there is no need to derive problem-specific weak formulations or to assemble and store global sparse matrices. More importantly, for nonlinear phase-field systems, strong-form methods avoid repeated Jacobian assembly and factorization during Newton iterations, which often becomes a dominant bottleneck for convergence, memory usage, and parallel computational complexity in weak-form FEM solvers. Instead, strong-form stencil methods apply local operators, reducing memory movement and simplifying data layout. However, despite their advantages in computational and data-storage complexity, existing strong-form approaches remain limited in their ability to robustly maintain formal accuracy on adaptive, non-uniform meshes for phase field modeling.

Phase-field dynamics are inherently interface-dominated and often involve curvature-sensitive or even fourth-order operators. Under these conditions, localized discretization defects translate directly into incorrect interface motion, spurious currents, or degraded convergence. FVM (finite volume method) is naturally flux-based for second-order operators, so fourth-order terms typically require mixed formulations or wide stencils, complicating boundary conditions, increasing communication/ghost layers on AMR, and making it difficult to achieve simultaneous discrete mass conservation and free-energy dissipation. In practice, these difficulties often necessitate auxiliary-variable reformulations [27], specialized coarse–fine flux matching [28], or correction and mortar strategies [29]. GFDM (generalized finite-different method) can approximate higher derivatives *via* local polynomial fitting and least-squares stencils, but on AMR these stencils become irregular near interfaces, degrading accuracy and isotropy and potentially producing spurious oscillations or interface pinning; moreover, conservation and energy stability are not automatic and can be violated by adaptive remeshing and refinement level transfer, which may inject high-frequency noise and alter mass/energy.

To address these challenges, we introduce an alternative strong-form based numerical framework for phase field modeling by coupling the Gradient Smoothing Method (GSM) with a hierarchical adaptive triangle-structured mesh. GSM constructs gradients and higher-order operators through geometric smoothing and local averaging, which can explicitly conserve quantities in approximation, rather than face-based flux reconstruction or polynomial fitting, making the operators far less sensitive to stencil irregularity, mesh anisotropy, and refinement transitions. We propose a hierarchical adaptive mesh refinement (HAMR) strategy that enforces uniform fine resolution within the interfacial band and coarse in bulk, while confining mesh non-uniformity to a narrow, well-controlled transition layer. This design allows the GSM solver to retain overall second-order accuracy and desired convergence rate. Remeshing is local and hierarchical: refinement and coarsening decisions are automatically linked to the physical interface thickness, remeshing operations are restricted to a small set of neighbors around the moving interface, and parent–child connectivity is stored during refinement so that inverse coarsening can be performed by direct rollback/merge operations without costly global searches.

The remainder of this paper is organized as follows. Section 2 summarizes the governing phase-field models. Section 3 introduces the Gradient Smoothing Method, including the stencil construction, discrete operators, numerical stability criterion, and the fully vectorized implementation of GSM Laplace operator. Section 4 describes the hierarchical adaptive meshing algorithm, its refinement criteria tied to physical properties, and parametric studies of mesh quality. Section 5 presents numerical experiments that validate accuracy, demonstrate efficiency and computational complexity against uniform finite-difference baselines, and Section 6 concludes.

## 2. Phase Field Method

The Phase Field Method involves two classical governing equations: the Allen–Cahn (A–C) equation and the Cahn–Hilliard (C–H) equation [9, 30], which are respectively suited for non-conserved and conserved order parameters, respectively.

The Allen–Cahn equation describes the temporal evolution of a non-conserved order parameter $\phi(\mathbf{r}, t)$, which typically represents phase identity. Specifically, $\phi = 1$ denotes one pure phase and $\phi = -1$ denotes the other, while $-1 < \phi < 1$ represents the diffusive interface between the two phases, as shown in Figure 1(a).

The phase-field formulation is derived from the variational principle of minimizing the total free energy functional

$$F(\phi) = \int_{\Omega} \left( f(\phi) + \frac{\kappa}{2} |\nabla\phi|^2 \right) d\Omega \tag{1}$$

where $f(\phi)$ is a bulk free-energy density with a double-well potential and $\kappa$ is the gradient energy coefficient. One commonly used double-well potential, as plotted in Figure 1(b), is

$$f(\phi) = f_0 \frac{(\phi - 1)^2 (\phi + 1)^2}{4} \tag{2}$$

with $f_0$ being a constant controlling the magnitude of the potential.

With this form of free energy, the governing A-C equation becomes

$$\frac{\partial\phi}{\partial t} = -M_\phi \left( f'^{(\phi)} - \kappa \nabla^2 \phi \right) = -M_\phi [f_0(\phi^3 - \phi) - \kappa \nabla^2 \phi] \tag{3}$$

where $M_\phi$ is mobility. The A-C equation drives the system toward local free-energy minimization and is widely used to model interface motion driven by curvature, phase transformation, and order-disorder transitions. Although computationally efficient, the Allen-Cahn equation does not conserve the total phase volume, which limits its applicability in problems requiring mass conservation.

The C-H equation governs the evolution of a conserved order parameter, such as concentration $c(\mathbf{r}, t)$. It is expressed as

$$\frac{\partial c}{\partial t} = \nabla \cdot (M_c \nabla \mu) \tag{4}$$

where $M_c$ is mobility and $\mu$ is the chemical potential defined as

$$\mu = f'(c) - \kappa \nabla^2 c = f_0(c^3 - c) - \kappa \nabla^2 c \tag{5}$$

with the simple double-well potential in Eq. (2).

By construction, the C–H equation preserves the global integral of order parameter over the computational domain, making it suitable for modeling phase separation [31-33], spinodal decomposition [34, 35], and diffusion-controlled coarsening processes [36, 37]. Compared with the A–C equation, the C–H equation is fourth-order in space, which significantly increases numerical stiffness and places higher demands on spatial discretization accuracy and stability.

In this work, our primary objective is to demonstrate the computational efficiency and robustness of the proposed GSM solver based on an adaptive structured mesh for phase-field modeling. Therefore, we do not focus on a specific physical application but instead consider the general forms of the A–C and C–H equations with the simple double-well potential introduced above.

For simplicity, the mobilities are assumed to be constant and are set to $M_\phi = 1$ and $M_c = 1$. In addition, $f_0 = 1$ is adopted. The governing equations are then reduced to

$$\frac{\partial\phi}{\partial t} = -(\phi^3 - \phi) + \kappa \nabla^2 \phi \tag{6}$$

and

$$\frac{\partial c}{\partial t} = \nabla^2 (c^3 - c - \kappa \nabla^2 c) \tag{7}$$

for the A-C and C-H models, respectively.

The gradient energy coefficient $\kappa$ controls the interface thickness, and it is worthwhile to mention that the above simple phase-field models admit an analytical equilibrium solution [38] that describes the smooth diffuse interface:

$$c(x) = \tanh\left(\frac{d}{\sqrt{2\kappa}}\right) \tag{8}$$

where $d$ denotes the distance from the interface center. If we define the interface thickness as the region where $-0.95 < c < 0.95$ in numerical simulations, then, after some algebraic deductions, we get that the interface thickness, $\delta$, can be expressed in terms of the gradient-energy coefficient $\kappa$:

$$\delta = 2\sqrt{2\kappa}\,\mathrm{atanh}(0.95) \approx 5.18\sqrt{\kappa} \tag{9}$$

That is, the interface thickness scales as $\delta \sim \sqrt{\kappa}$, implying that the interface thickness doubles when $\kappa$ is increased by a factor of four. In this study, different interface thicknesses are generated by varying $\kappa$, allowing us to systematically evaluate the robustness and accuracy of the proposed adaptive mesh-based GSM solver under different interfacial resolutions. The relationship in Eq. (9) can be used to theoretically determine the mesh resolution required to properly resolve the interface for a given $\kappa$.

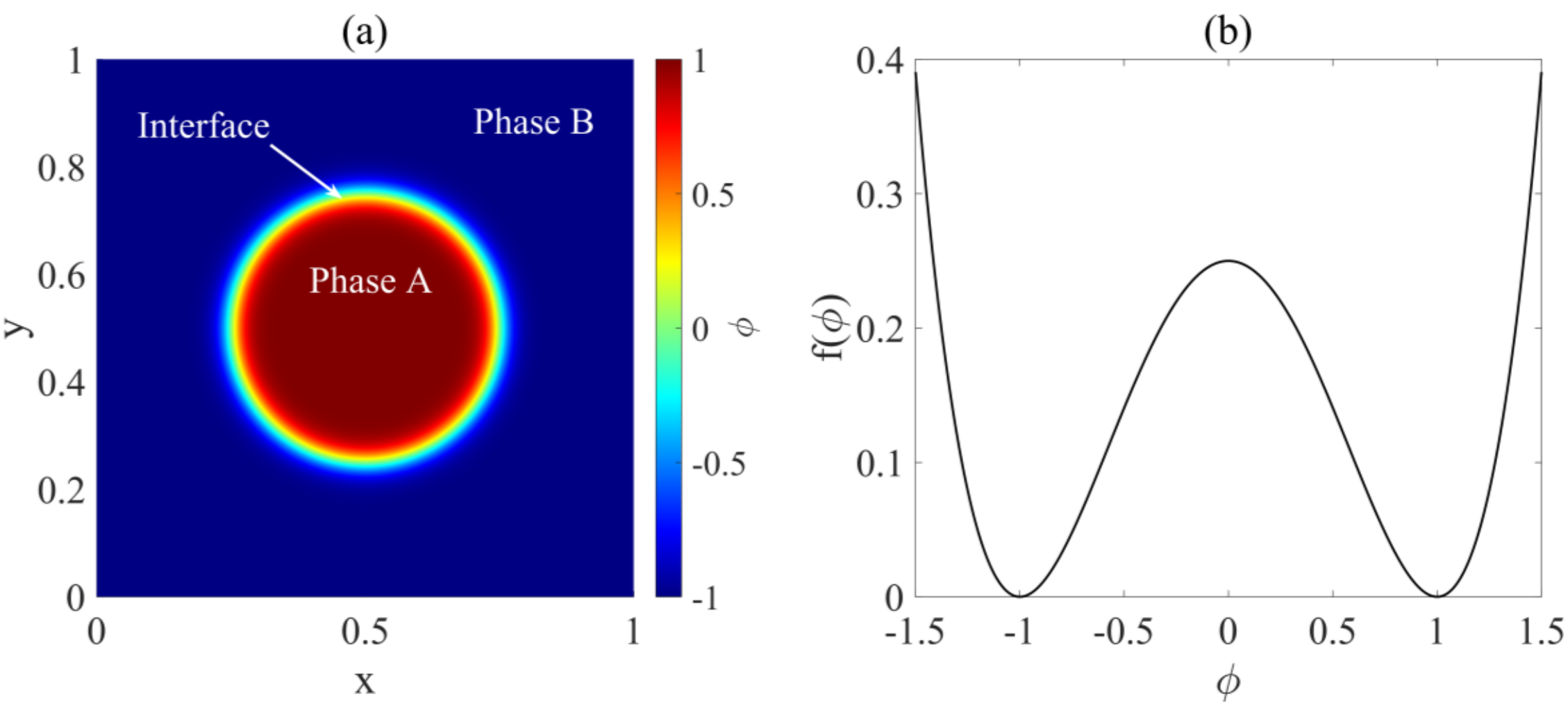

Figure 1 (a) Order-parameter field for a two-phase A/B system with a diffuse interface representation commonly used in phase-field formulations. (b) The double-well free-energy potential $f(\phi)$ with two minima representing the two equilibrium phase states and enforcing phase separation. $f_0 = 1$.

### 3. Gradient Smoothing theory

Consider a material domain $\Omega$ bounded by $\partial\Omega$. In the GSM framework, the gradient of a field variable $U(\mathbf{r})$ at any location $\mathbf{r}$ is approximated using a smoothing operator defined over a local smoothing domain $\Omega^S$ surrounding the location $\mathbf{r}$, as shown in Figure 2(a). The GSM gradient operator is defined by the integral representation

$$\nabla U(\mathbf{r}) = \int_{\Omega^S} \nabla U(\mathbf{r}')\widehat{w}(\mathbf{r}-\mathbf{r}')\mathrm{d}\mathbf{r}' \tag{10}$$

where $\Omega^S$ is referred as gradient smoothing domain or GSD, and $\widehat{w}$ is a smoothing function that must be sufficiently smooth and at least first order differentiable within the GSD. The accuracy of Eq. (10) depends on the choice of the smoothing function $\widehat{w}$ [39]. A fundamental requirement is that a constant gradient field be exactly reproduced, leading to the unity condition

$$\int_{\Omega^S} \widehat{w}(\mathbf{r}-\mathbf{r}')\mathrm{d}\mathbf{r}' = 1 \tag{11}$$

In this work, the simplest constant smoothing function is adopted:

$$\widehat{w} = \begin{cases} 0, & \mathbf{r} \notin \Omega^S \\ \dfrac{1}{V}, & \mathbf{r} \in \Omega^S \end{cases} \tag{12}$$

where $V$ denotes the area of $\Omega^S$.

Applying integration by parts to Eq. (10) yields

$$\nabla U(\mathbf{r}) = \int_{\partial\Omega^S} U(\mathbf{r}')\widehat{w}(\mathbf{r}-\mathbf{r}')\vec{n}\mathrm{d}S - \int_{\Omega^S} U(\mathbf{r}')\nabla\widehat{w}(\mathbf{r}-\mathbf{r}')\mathrm{d}\mathbf{r}' \tag{13}$$

where $\partial\Omega^S$ is its boundary, as shown in Figure 2(a), and $\vec{n}$ is the outwards unit normal vector on $\partial\Omega^S$. Since $\widehat{w}$ is constant within $\Omega^S$, its gradient vanishes, and Eq. (13) simplifies to

$$\nabla U(\mathbf{r}) = \frac{1}{V_i}\int_{\partial\Omega^S} U(\mathbf{r}')\vec{n}\mathrm{d}S \tag{14}$$

To evaluate the boundary integral in Eq. (14), different discretization schemes have been proposed in the literature [24]. In this study, the one-point quadrature (rectangular rule), based on the node gradient smoothing domain (nGSD) as illustrated in Figure 2(b), is adopted due to its simplicity and sufficient accuracy.

Using the rectangular rule, the boundary integral is approximated by summation over the boundary segments of the nGSD. The field value $U$ on each segment is evaluated at the midpoint between the target particle $i$ and its supporting neighbor $j_k$. Consequently, the GSM gradient at particle $i$ is discretized as

$$\begin{cases} \dfrac{\partial U_i}{\partial x} \approx \dfrac{1}{V_i}\displaystyle\sum_{k=1}^{n_i} U_{m_k}(\Delta S_x)_{ij_k} \\ \dfrac{\partial U_i}{\partial y} \approx \dfrac{1}{V_i}\displaystyle\sum_{k=1}^{n_i} U_{m_k}(\Delta S_y)_{ij_k} \end{cases} \tag{15}$$

where $V_i$ is the area of the nGSD for particle $i$, and $n_i$ is the number of supporting neighbors. $U_{m_k}$ represents the field variable $U$ at the midpoint $m_k$ between the particles $i$ and $j_k$, which is evaluated by linear interpolation:

$$U_{m_k} = \frac{U_i + U_{j_k}}{2} \tag{16}$$

The shape function components $\Delta(S_x)_{ij_k}$ and $\Delta(S_y)_{ij_k}$ are given by

$$\Delta(S_x)_{ij_k} = \Delta(S_x)_{ij_k}^{(L)} + \Delta(S_x)_{ij_k}^{(R)}, \qquad \Delta(S_y)_{ij_k} = \Delta(S_y)_{ij_k}^{(L)} + \Delta(S_y)_{ij_k}^{(R)} \tag{17}$$

with

$$\begin{cases} \Delta(S_x)^{(L)}_{ij_k} = l_{m_k c_k}(n_x)_{m_k c_k}, & \Delta(S_y)^{(L)}_{ij_k} = l_{m_k c_k}(n_y)_{m_k c_k} \\ \Delta(S_x)^{(R)}_{ij_k} = l_{m_k c_{k-1}}(n_x)_{m_k c_{k-1}}, & \Delta(S_y)^{(R)}_{ij_k} = l_{m_k c_{k-1}}(n_y)_{m_k c_{k-1}} \end{cases} \tag{18}$$

where $l$ is edge length, $(n_x, n_y)$ are the components of the outward unit vector $\vec{n}$, and $c_k$ is the centroid of the triangle formed by points $i$, $j_k$, and $j_{k+1}$, as shown in Figure 2(b).

Equations (15)-(18) define the GSM gradient smoothing operator for a field variable $U$ at particle $i$. By successively applying the gradient smoothing operator twice, we get the GSM Laplace operator

$$\nabla \cdot (\nabla U_i) \approx \nabla \cdot (\nabla U_i)|_{\text{GSM}} = \frac{1}{V_i} \sum_{k=1}^{n_i} \left[ \frac{\partial U_{m_k}}{\partial x} (\Delta S_x)_{ij_k} + \frac{\partial U_{m_k}}{\partial y} (\Delta S_y)_{ij_k} \right] \tag{19}$$

Note that, a linear interpolation of midpoint gradient

$$\frac{\partial U_{m_k}}{\partial x} = \frac{\frac{\partial U_i}{\partial x} + \frac{\partial U_j}{\partial x}}{2} \quad and \quad \frac{\partial U_{m_k}}{\partial y} = \frac{\frac{\partial U_i}{\partial y} + \frac{\partial U_j}{\partial y}}{2} \tag{20}$$

will lead to the so-called checkerboard instability [40], which destroys the accuracy of the Laplace operator. To eliminate this issue, a directional correction [40] is necessarily applied:

$$\nabla U_{m_k} = \nabla U_{m_k} - \left[ \nabla U_{m_k} \cdot \vec{t}_{ij} - \frac{U_j - U_i}{d_{ij}} \right] \vec{t}_{ij} \tag{21}$$

where $\vec{t}_{ij}$ is the unit vector along the line connecting particles $i$ and $j$, and $d_{ij}$ is their distance.

The GSM Laplace operator in Eq. (19), together with the directional correction in Eq. (21), is employed to discretize the second-order derivative terms in the Allen–Cahn and Cahn–Hilliard equations. This formulation performs robustly on irregular meshes. Theoretical analyses have indicated that this form of GSM Laplace operator achieves second-order accuracy on uniform meshes and first-order accuracy on non-uniform meshes [24, 26], which will be verified numerically in subsequent section 5.

Note that, considering the formulations of the geometric edge coefficients in Eqs. (17) and (18), the flux contribution of the neighboring particle $j$ to particle $i$ is always antisymmetric to that of $i$ to $j$, *i.e.*, $\mu_{m_k}(\Delta S_x)_{ij_k} \equiv -\mu_{m_k}(\Delta S_x)_{j_k i}$ because $(\Delta S_x)_{ij_k} \equiv -(\Delta S_x)_{j_k i}$. This guarantees exact conservation of the chemical potential flux in the C-H model, and enforces strict concentration conservation in GSM simulations.

When explicit time integration is used, the stability constraints for the GSM solver follow those of finite difference schemes. In two-dimensional simulations, the time step must satisfy

$$\Delta t_{AC} \leq \frac{h^2}{4\kappa} \tag{22}$$

for the A-C equation, and

$$\Delta t_{CH} \leq \frac{h^2}{4 + \frac{8\kappa \times 4}{h^2}} \tag{23}$$

for the C-H equation [41, 42], where $h$ denotes the grid size.

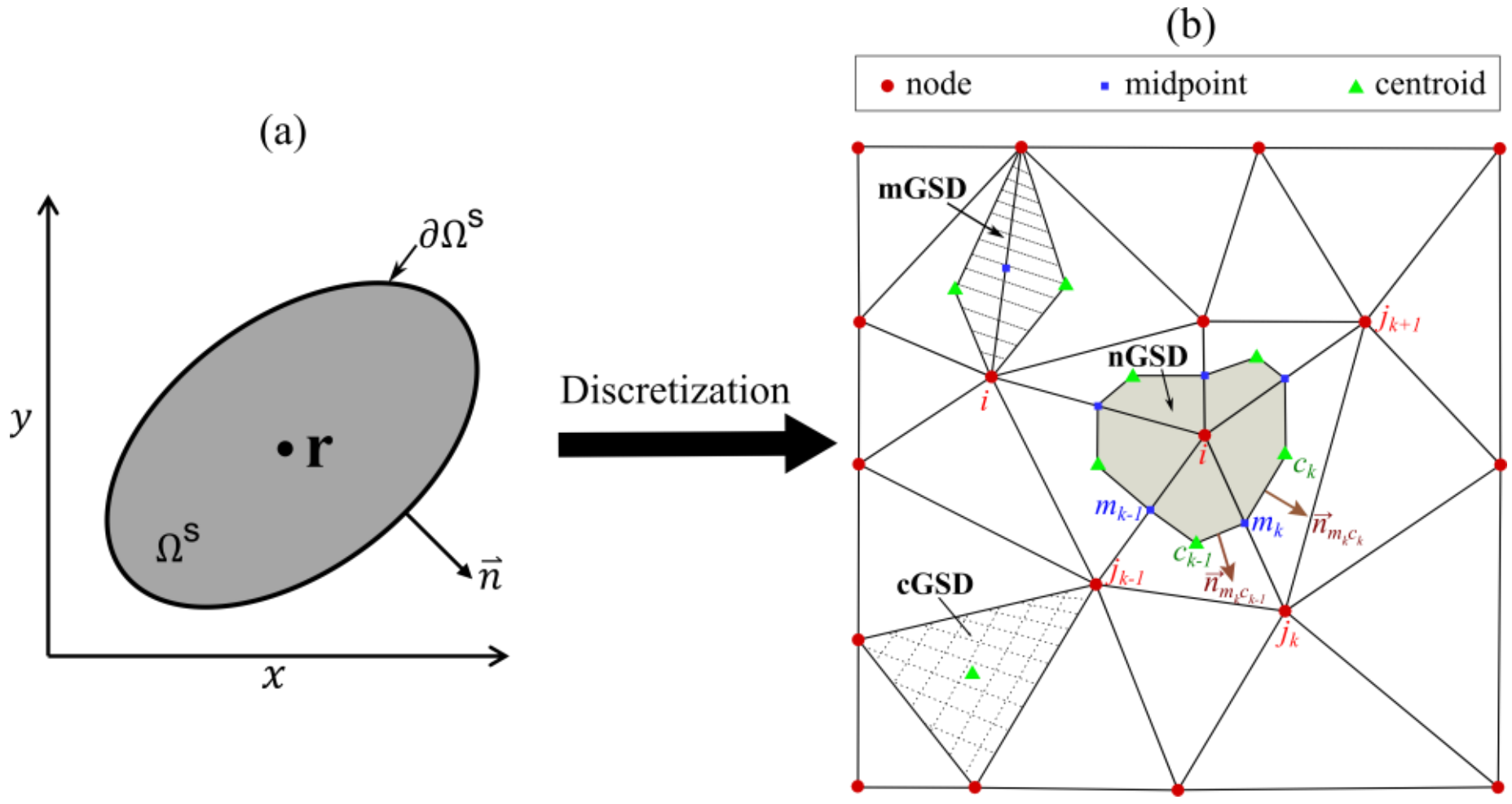

Figure 2 Illustration of generic gradient smoothing domain (a) and the discretized smoothing domains (b) in GSM. nGSD, mGSD, and cGSD denote the gradient smoothing domain (GSD) of node, midpoint, and centroid, respectively.

*Full vectorization of the GSM Laplacian on an adaptive triangular mesh*

To optimize computational efficiency, we implement a fully vectorized GSM-based Laplacian operator on an adaptively refined, structured triangular mesh.

The domain is discretized by a triangulation with nodal coordinates arrays

$$\mathbf{x} \in \mathbb{R}^{N_d}, \mathbf{y} \in \mathbb{R}^{N_d} \tag{24}$$

and an element connectivity array

$$\mathbf{T} \in \mathbb{N}^{N_T \times 3}, \qquad \mathbf{T}(e,:) = (v_{e1}, v_{e2}, v_{e3}) \tag{25}$$

where $N_d$ and $N_T$ denote the numbers of nodes and triangular elements, respectively. $(v_{e1}, v_{e2}, v_{e3})$ are the global vertex indices of triangle element $e$.

For a nodal scalar field $\mathbf{u} \in \mathbb{R}^{N_d}$ associated with the nodal coordinates $\mathbf{x}$ and $\mathbf{y}$, the objective is to compute a nodal approximation to $\Delta\mathbf{u}$ without element-by-element loops, by evaluating all per-triangle geometric quantities in matrix form and assembling nodal quantities *via* vectorized scatter-add operations.

*Triangle-wise gather.* Using the connectivity $\mathbf{T}$, nodal coordinates and field variables are gathered into triangular-wise matrices

$$\mathbf{X} = \mathbf{x}(\mathbf{T}) \in \mathbb{R}^{N_T \times 3}, \qquad \mathbf{Y} = \mathbf{y}(\mathbf{T}) \in \mathbb{R}^{N_T \times 3}, \qquad \mathbf{U} = \mathbf{u}(\mathbf{T}) \in \mathbb{R}^{N_T \times 3} \tag{26}$$

whose columns correspond to the three local vertices of each triangle. This gathering step produces length-$N_T$ arrays for each local quantity and enables fully vectorized evaluation of edge-based coefficients.

*Mid-edge values.* The GSM flux construction uses mid-edge averages (computed triangle-wise for all triangles) on the three edges $\{12, 23, 31\}$:

$$\mathbf{U}_{12} = \frac{\mathbf{U}_1 + \mathbf{U}_2}{2}, \qquad \mathbf{U}_{23} = \frac{\mathbf{U}_2 + \mathbf{U}_3}{2}, \qquad \mathbf{U}_{31} = \frac{\mathbf{U}_3 + \mathbf{U}_1}{2} \tag{27}$$

with

$$\mathbf{U}_1 = \mathbf{U}(:,1), \qquad \mathbf{U}_2 = \mathbf{U}(:,2), \qquad \mathbf{U}_3 = \mathbf{U}(:,3) \tag{28}$$

*GSM geometric edge coefficients*. For each triangle, we define edge coefficients ($\Delta(S_x)$ and $\Delta(S_y)$) on the three edges $\{12, 23, 31\}$:

$$\begin{cases} \Delta(\mathbf{S}_x)_{12} = \dfrac{2\mathbf{Y}_3 - \mathbf{Y}_1 - \mathbf{Y}_2}{6}, & \Delta(\mathbf{S}_y)_{12} = \dfrac{2\mathbf{X}_3 - \mathbf{X}_1 - \mathbf{X}_2}{6} \\ \Delta(\mathbf{S}_x)_{23} = \dfrac{2\mathbf{Y}_1 - \mathbf{Y}_2 - \mathbf{Y}_3}{6}, & \Delta(\mathbf{S}_y)_{23} = \dfrac{2\mathbf{X}_1 - \mathbf{X}_2 - \mathbf{X}_3}{6} \\ \Delta(\mathbf{S}_x)_{31} = \dfrac{2\mathbf{Y}_2 - \mathbf{Y}_3 - \mathbf{Y}_1}{6}, & \Delta(\mathbf{S}_y)_{31} = \dfrac{2\mathbf{X}_2 - \mathbf{X}_3 - \mathbf{X}_1}{6} \end{cases} \tag{29}$$

with

$$\begin{cases} \mathbf{X}_1 = \mathbf{X}(:,1), & \mathbf{X}_2 = \mathbf{X}(:,2), & \mathbf{X}_3 = \mathbf{X}(:,3) \\ \mathbf{Y}_1 = \mathbf{Y}(:,1), & \mathbf{Y}_2 = \mathbf{Y}(:,2), & \mathbf{Y}_3 = \mathbf{Y}(:,3) \end{cases} \tag{30}$$

These coefficients act as precomputed geometric weights that convert edge-centered quantities into vertex-centered divergence-like balances.

*Normalization factor matrix* **V**. In the vectorized code, the nodal normalization factor (corresponding to the area of GSD, $V$, in Eq. (19)) is assembled by accumulating a per-triangle weight

$$V_e = \frac{h_c^2}{6 \times 2^{l_e}} \tag{31}$$

to each of the triangle's three vertices, where $h_c$ is the global grid size of coarsest mesh, and $l_e$ is the refinement level of triangle element $e$. The resultant nodal scale for node $i$ is

$$V_i = \sum_{e:\, i \in e} V_e \tag{32}$$

where $V_i$ reflects the area of node $i$ and is the $i$-th element of the normalization factor matrix **v**.

*Vectorized GSM gradient assembly*. Using mid-edge values and the geometric coefficients, we form per-edge flux contribution as:

$$\mathbf{C}_{ab}^x = \mathbf{U}_{ab} @ \Delta(\mathbf{S}_x)_{ab}, \;\; \mathbf{C}_{ab}^y = \mathbf{U}_{ab} @ \Delta(\mathbf{S}_y)_{ab}, \quad ab \in \{12, 23, 31\} \tag{33}$$

where @ denotes element-wise multiplication.

These terms are assembled to the three vertices of each triangle using the fixed GSM divergence pattern employed in the code:

$$\begin{cases} \mathbf{G}_1^x = \mathbf{C}_{12}^x - \mathbf{C}_{31}^x, & \mathbf{G}_1^y = \mathbf{C}_{12}^y - \mathbf{C}_{31}^y \\ \mathbf{G}_2^x = \mathbf{C}_{23}^x - \mathbf{C}_{12}^x, & \mathbf{G}_2^y = \mathbf{C}_{23}^y - \mathbf{C}_{12}^y \\ \mathbf{G}_3^x = \mathbf{C}_{31}^x - \mathbf{C}_{23}^x, & \mathbf{G}_3^y = \mathbf{C}_{31}^y - \mathbf{C}_{23}^y \end{cases} \tag{34}$$

Followed by a scatter-add over all triangles to obtain the global nodal accumulations ($\mathbf{g}_x, \mathbf{g}_y$)

$$\begin{cases} \mathbf{g}_x(\mathbf{v}_1) += \mathbf{G}_1^x, & \mathbf{g}_y(\mathbf{v}_1) += \mathbf{G}_1^y \\ \mathbf{g}_x(\mathbf{v}_2) += \mathbf{G}_2^x, & \mathbf{g}_y(\mathbf{v}_2) += \mathbf{G}_2^y \\ \mathbf{g}_x(\mathbf{v}_3) += \mathbf{G}_3^x, & \mathbf{g}_y(\mathbf{v}_3) += \mathbf{G}_3^y \end{cases} \tag{35}$$

recall that $(\mathbf{v}_1, \mathbf{v}_2, \mathbf{v}_3)$ are the three columns of connectivity matrix $\mathbf{T}$.

The nodal gradient is then computed by normalized and eventually summed up to the corresponding nodes:

$$\frac{\partial \mathbf{u}}{\partial x} = \frac{\mathbf{g}_x}{\mathbf{v}}, \qquad \frac{\partial \mathbf{u}}{\partial y} = \frac{\mathbf{g}_y}{\mathbf{v}} \tag{36}$$

based on Eq. (15).

To build the GSM Laplacian operator, the same GSM assembly is applied again, but with the scalar midpoint values $\mathbf{U}_{ab}$ replaced by midpoint gradients. In implementation, midpoint gradients are first obtained by averaging nodal gradients along each edge:

$$\frac{\partial \mathbf{U}_{ab}}{\partial x} = \frac{\frac{\partial \mathbf{U}_a}{\partial x} + \frac{\partial \mathbf{U}_b}{\partial x}}{2}, \qquad \frac{\partial \mathbf{U}_{m_k}}{\partial y} = \frac{\frac{\partial \mathbf{U}_a}{\partial y} + \frac{\partial \mathbf{U}_b}{\partial y}}{2} \tag{37}$$

To obtain an accurate GSM Laplace operator, the directional correction in Eq. (21) is applied to each edge midpoint gradient. Specifically, for edge $ab$ in each triangle, we define its length $\mathbf{l}_{ab}$ and unit tangent $\mathbf{t}_{ab}$ as

$$\mathbf{l}_{ab} = \sqrt{(\mathbf{X}_a - \mathbf{X}_b)^2 + (\mathbf{Y}_a - \mathbf{Y}_b)^2}, \qquad \mathbf{t}_{ab} = \frac{1}{\mathbf{l}_{ab}} @ \begin{bmatrix} \mathbf{X}_b - \mathbf{X}_a \\ \mathbf{Y}_b - \mathbf{Y}_a \end{bmatrix}, \qquad ab \in \{12, 23, 31\} \tag{38}$$

Let the uncorrected midpoint gradient be

$$\mathbf{U}'_{ab} = \begin{bmatrix} \frac{\partial \mathbf{U}_{ab}}{\partial x} \\ \frac{\partial \mathbf{U}_{ab}}{\partial y} \end{bmatrix} \tag{39}$$

The corrected GSM midpoint gradient used in Laplacian assembly is then

$$\mathbf{U}'^{c}_{ab} = \begin{bmatrix} \mathbf{U}^{x,c}_{ab} \\ \mathbf{U}^{y,c}_{ab} \end{bmatrix} = \mathbf{U}'_{ab} - \left[ \mathbf{U}'_{ab} \cdot \mathbf{t}_{ab} - \frac{\mathbf{U}_b - \mathbf{U}_a}{\mathbf{l}_{ab}} \right] \mathbf{t}_{ab} \tag{40}$$

Finally, the Laplacian term is assembled by applying the same GSM divergence pattern in Eq. (33) to the corrected edge gradients. Corrected flux components are defined as

$$\mathbf{C}^{x,c}_{ab} = \mathbf{U}^{x,c}_{ab} @ \Delta(\mathbf{S}_x)_{ab}, \ \mathbf{C}^{y,c}_{ab} = \mathbf{U}^{x,c}_{ab} @ \Delta(\mathbf{S}_y)_{ab}, \quad ab \in \{12, 23, 31\} \tag{41}$$

by replacing the field variable matrix $\mathbf{U}_{ab}$ in Eq. (33) with the corrected midpoint gradients. Following the derivations in Eqs. (34)-(35), we can get the nodal flux $\mathbf{g}_x^c$ and $\mathbf{g}_y^c$, and the GSM Laplace operator becomes:

$$\Delta \mathbf{u} = \frac{\partial^2 \mathbf{u}}{\partial x^2} + \frac{\partial^2 \mathbf{u}}{\partial x^2} = \frac{\mathbf{g}_x^c}{\mathbf{v}} + \frac{\mathbf{g}_y^c}{\mathbf{v}} \tag{42}$$

In terms of computational complexity, the vectorized implementation is dominated by triangle-wise array operations (gathers and per-triangle algebra on arrays with a dimension size of $N_T$) and a small number of global scatter-add assemblies (to build $\mathbf{v}$ in Eq. (32), and to assemble the gradient and corrected-divergence terms in Eq. (35)). This eliminates explicit element loops and yields an efficient Laplacian evaluation suitable for adaptive meshes.

## 4. **Hierarchical Adaptive Mesh Refinement (HAMR) algorithm**

This HAMR algorithm constructs an adaptive structured triangular mesh that automatically tracks an evolving interface described by a scalar indicator field, *i.e.*, the phase-field variable $\phi(\mathbf{r}, t)$ or $c(\mathbf{r}, t)$. The mesh is adaptively refined in the interfacial zone and coarsened in the bulk while remaining conforming, *i.e.*, free of hanging nodes. The conforming requirement is crucial for GSM solvers because it calculates smoothed derivatives from the flux/gradient information integrated over element (or smoothing-domain) boundaries, which implicitly assumes that adjacent elements share the same geometric edge and have compatible traces of the discrete field along that edge. Hanging nodes create non-matching coarse–fine interfaces, so the two sides of a nominally shared edge do not represent the same interpolated field, degrading $C_0$ continuity and corrupting the boundary data used in smoothing. This leads to biased smoothed derivatives near refinement transitions and can reduce accuracy and convergence conditions of GSM solvers. A triangular mesh is therefore adopted, together with a refinement strategy that explicitly enforces conformity. The adaptive triangular mesh is composed of different sizes of regular right isosceles triangles aligned with a Cartesian lattice, which avoids complicated global remeshing (*e.g.*, Delaunay retriangulation). Finally, a hierarchical representation (parent–child-neighbor relations and refinement levels) enables efficient local refinement/coarsening and neighbor queries by traversing the refinement tree rather than scanning the full mesh for searching neighboring elements, which is advantageous for efficiency and future GPU parallel implementations.

### *4.1 Intial mesh and interface indicator*

The computational domain is initialized with a coarse, uniform structured triangular mesh, as shown in Figure 3(a). The coarsest mesh resolution is controlled by $N_0$, the number of nodal points along each coordinate direction, and therefore $N_0$ sets the largest admissible element size in the domain. Another key parameter is the maximum refinement level, $\bar{n}_{level}$, which prescribes the minimum element size and the maximum refinement depth of the adaptive mesh, as shown in Figure 4. The adaptive strategy is designed to concentrate resolution only where it is needed. Specifically, triangles are refined only in the interfacial zone, where the evolution dynamics are significant, while the bulk region without dynamics is kept coarse. Thus, the adaptive mesh can be viewed as being bounded by two limiting uniform meshes: the coarsest mesh determined by $N_0$ (Figure 3a) and the finest mesh implied by $\bar{n}_{level}$ (Figure 3b and Figure 3c). At any time, the actual mesh is obtained by locally activating refinement between these limits according to the refinement indicator, producing a gradual transition from the coarse bulk mesh to the fine interfacial mesh.

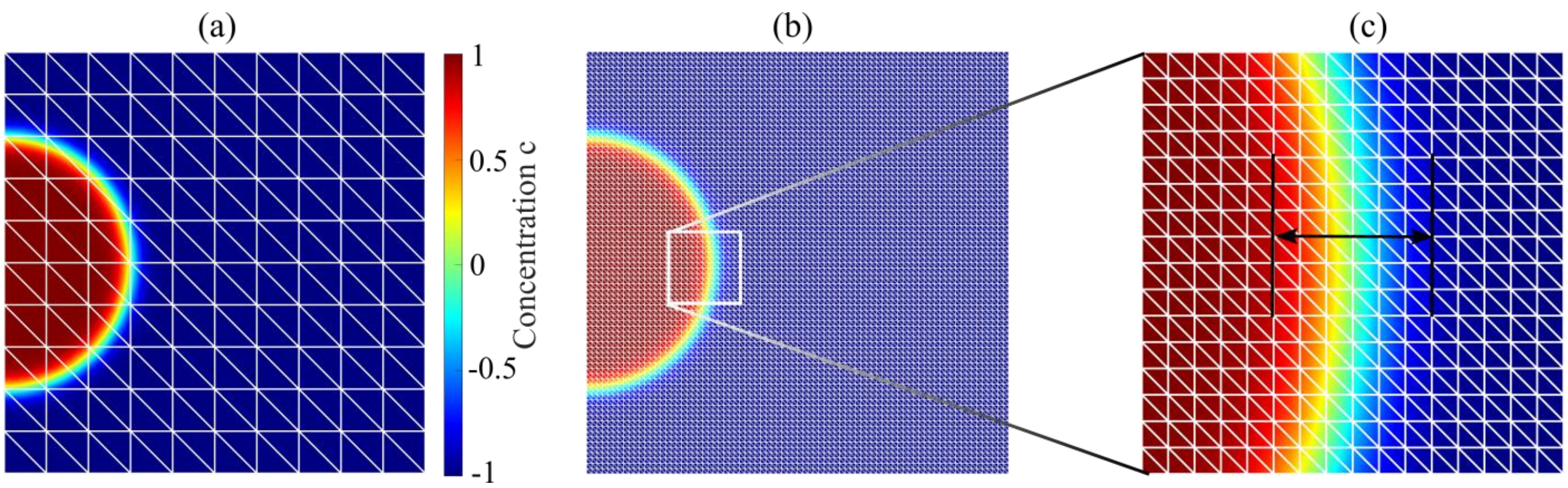

Figure 3 Coarsest and finest uniform meshes in adaptive structured mesh. (a) Coarsest uniform mesh controlled by $N_0$, which sets the largest element size in the bulk. (b) Finest uniform mesh controlled by $\bar{n}_{level}$, which defines the minimum element size near the interface. (c) Zoom-in of panel (b) showing the finest mesh near the diffuse interface in the concentration field.

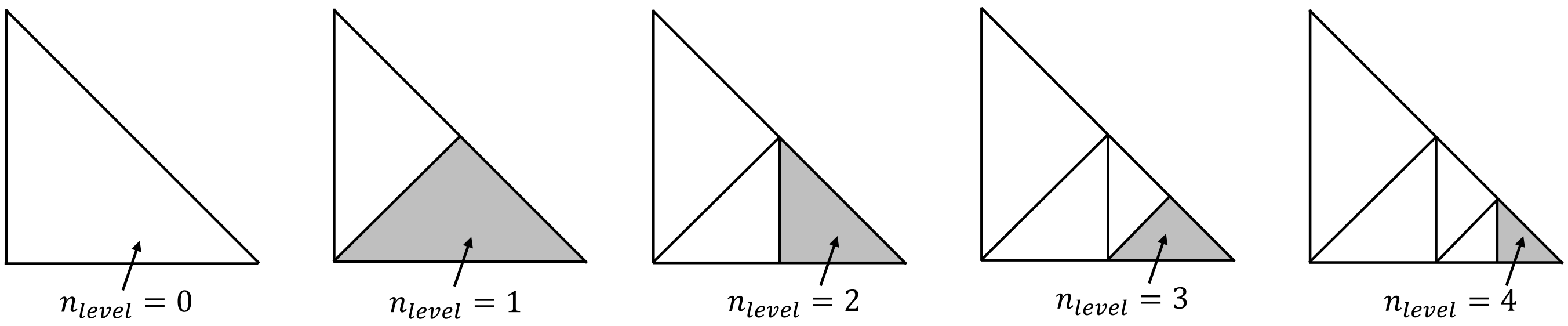

Figure 4 Refinement level $n_{level}$ of triangular elements employed in HAMR algorithm.

### *4.2 Hierarchical meshing*

The adaptive triangulation is represented with a hierarchical, tree-based data structure that enables local refinement and coarsening without global remeshing and searching.

In implementation, element-wise information is assembled into a matrix denoted by **T**, where each row stores all the necessary information corresponding to a triangular element. Specifically, the stored information includes: (i) the indices of its three vertices; (ii) the index of its mother element (set to zero for elements on the coarsest level), as illustrated in Figure 5(a); (iii) the index of its sister element generated from the same bisection (Figure 5a); (iv) its child status and child indices (two children if the element is refined, as shown in Figure 5a); (v) its grandchild status (zero if none, one if present), as illustrated in Figure 5(b); (vi) the indices of its three edge neighbors (Figure 5c), with a value of zero indicating a boundary edge; and (vii) its refinement level. As illustrated in Figure 5, mother triangle 1 is bisected along its long edge to generate child triangles 5 and 6. If triangle 1 belongs to the coarsest mesh level, then triangles 1–3 have a refinement level of zero, indicating no refinement. The triangles 5 and 6 have a refinement level of 1, while the triangles 7 and 8 are grandchildren of triangle 1 and have refinement level 2. Thus, for any refinement level, a grandchild triangle has half the grid size and a refinement level two higher than its grandmother.

Nodal data consists of vertex coordinates **r** and the associated field variables **u** (*e.g.,* the phase-field order parameter). During refinement, new mid-edge nodes are appended to the nodal list, while coarsening removes nodes that are no longer required once child elements are merged back into their parent.

This hierarchical structure provides a persistent mapping between coarse- and fine-scale elements and enables coarsening as a controlled reversal of refinement, substantially reducing computational cost.

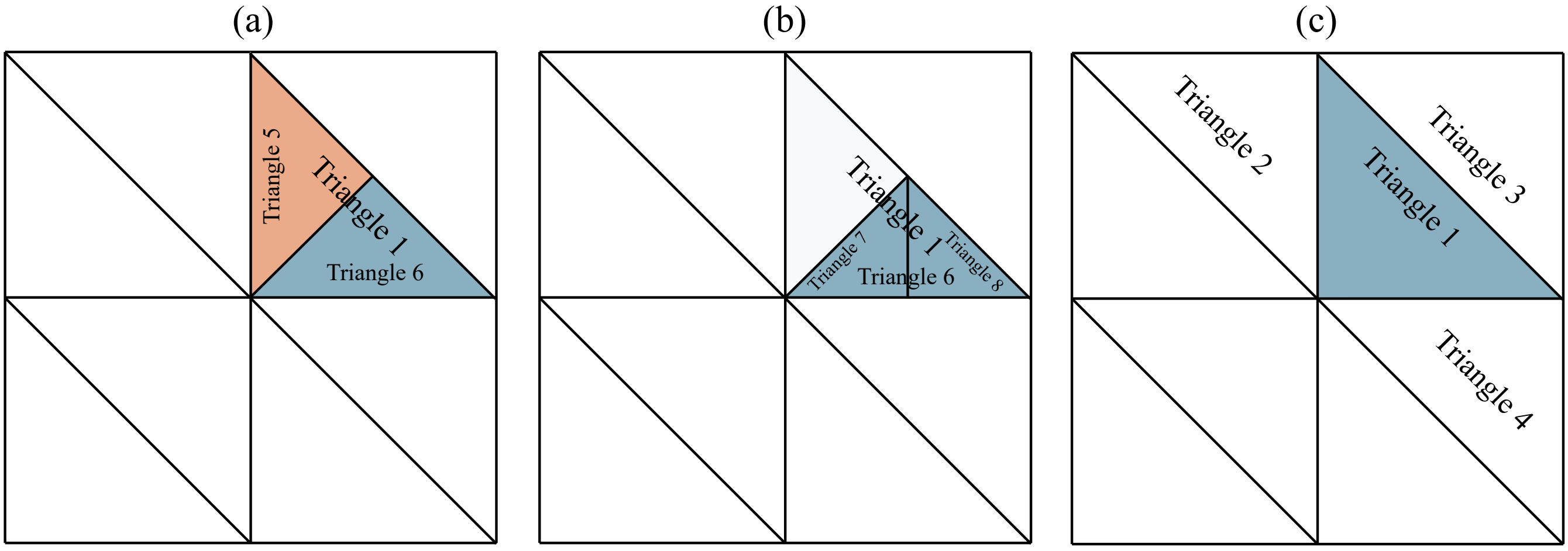

Figure 5 Illustration of hierarchical triangular element relationships in a structured right-isosceles mesh. (a) Mother triangle and its two child triangles generated by longest-edge bisection; triangles 5 and 6 are sister elements created from refining triangle 1. (b) Grandchildren elements; triangles 7 and 8 are obtained by refining child elements of triangle 1. (c) Edge-neighbor relationships; triangle 1 has three edge neighbors (triangles 2, 3, and 4), each sharing one edge with triangle 1.

### *4.3 Refinement strategy*

The refinement strategy is designed to accurately resolve evolving interfaces while maintaining mesh conformity, locality, and computational efficiency. Since phase-field interfacial problems typically exhibit sharp gradients confined to narrow regions, refinement is restricted to elements where solution variation exceeds a prescribed threshold or identified interfacial region. The refinement strategy is governed by three core principles: (i) deterministic refinement geometry, ensured by consistent vertex ordering and longest-edge selection; (ii) solution-driven refinement indicator, based on element-wise field variables; (iii) local conformity enforcement, achieved through refinement propagation across principal neighbors. Each refinement operation splits a mother triangle into two child elements, forming a refinement tree in which every active element explicitly records its ancestry and descendants.

*Deterministic longest-edge identification.* To ensure unique and repeatable refinement patterns, the vertices of each triangle are ordered in a counter-clockwise manner, with the right-angle vertex placed first. Let the ordered vertex indices of element $K$ be $(v_1, v_2, v_3)$. The longest edge is therefore

$$e_{\max}(K) = (v_2, v_3) \tag{43}$$

Under this convention, the second and third vertices always define the longest edge of the triangle. The neighboring element sharing this edge is uniquely identified as its *principal neighbor*. This deterministic definition eliminates ambiguity in refinement propagation and guarantees mesh conformity without introducing hanging nodes.

*Element-wise refinement indicator.* Let $u(\mathbf{r})$ denote the scalar field of interest (*e.g.*, phase-field order parameters $\phi$ and $c$), whose values are stored in the phase-field matrix $\mathbf{u}$, *i.e.*, $u_i = u(\mathbf{r}_i)$. For a triangular element $K$ with vertices $\{i, j, k\}$, an element-wise refinement indicator is defined as

$$\eta_K = \frac{u_i + u_j + u_k}{3} \tag{44}$$

An element is *marked* for refinement if

$$|\eta_K| < u_{thr} \tag{45}$$

where $u_{thr}$ is a user-defined tolerance defining interface region. This criterion identifies interfacial regions undergoing strong solution evolution. In particular, it drives all elements intersecting the interfacial region to the prescribed maximum refinement level $\bar{n}_{level}$, resulting in a locally uniform fine mesh. Such uniform refinement near the interface is essential for optimizing the numerical accuracy of the GSM solver, which relies on consistent element resolution to achieve high order accuracy. For this reason, we adopt this refinement indicator rather than the commonly used order-parameter–variation criterion, which does not ensure a uniformly fine mesh at the interface and can therefore degrade the GSM solver's high-order accuracy.

For each *marked* element $K$, refinement is performed by bisecting its longest edge. Let $\mathbf{r}_m = (x_m, \ y_m)$ denotes the midpoint of the longest edge,

$$x_m = \frac{1}{2}\left(x_{v_2} + x_{v_3}\right), \qquad y_m = \frac{1}{2}\left(y_{v_2} + y_{v_3}\right) \tag{46}$$

The mother element $K$ is split into two child elements $K_1$ and $K_2$, each inheriting one half of the longest edge and sharing the newly created midpoint. Both children are assigned a refinement level increased by one relative to the mother element. When updating the nodal coordinate matrix $\mathbf{r}$, the phase-field variable matrix $\mathbf{u}$ is updated accordingly by appending the order parameter value at the new midpoint through linear interpolation:

$$u_m = \frac{1}{2}\left(u_{v_2} + u_{v_3}\right) \tag{47}$$

with $u_m$ being the field variable $u$ at the newly appended midpoint.

This operation updates the hierarchical data structure by: (i) recording child indices and sister relationships, (ii) marking the mother element as inactive, and (iii) updating ancestry and refinement level information.

*Conformity enforcement via refinement propagation.* To preserve mesh conformity, refinement may need to propagate to neighboring elements. After bisecting a marked element, its *principal neighbor* is examined. If the principal neighbor is coarser and unrefined, it is marked for refinement to prevent hanging nodes. This process is repeated iteratively until all refined elements satisfy the conformity condition. Due to the hierarchical storage and explicit neighbor indexing, refinement propagation is confined to a small, localized set of elements.

The complete refinement procedure is summarized as Algorithm 1 in Figure 6.

```
Algorithm 1: Mesh refinement via longest-edge bisection
  Input: Element matrix T, node coordinates r, nodal values u
  Output: Updated T, r, u
  Compute element-wise indicators η_K for all active elements;
  Mark elements: R ← (|η_K| < u_thr) ∧ (T.child_flag = 0);
  while any(R) do
  |  /* Longest-edge refinement                                   */
  |  Identify longest edge and principal neighbor for elements in R;
  |  Create mid-edge nodes and append to nodal list;
  |  r ← append_nodes(r, mid-edge nodes);
  |  u ← append_values_by_interpolation(u, r, mid-edge nodes);
  |  Split each marked mother element into two children;
  |  Update child indices, sister links, and refinement levels;
  |  Mark mother elements as inactive;
  |  /* Conformity enforcement                                     */
  |  Identify principal neighbors of newly refined elements;
  |_ R ← (principal neighbor is active, unrefined, and coarser);
  Update hierarchy flags (children/grandchildren/sister/neighbor);
  Update local neighbor connectivity;
```

Figure 6 Pseudocode for hierarchical adaptive mesh refinement using longest-edge bisection.

### *4.4 Coarsening strategy*

Coarsening is implemented as a conservative reversal of refinement, enabled by the hierarchical tree-based mesh representation. This substantially reduces computational cost while maintaining conformity, consistency, and numerical accuracy. The operation is strictly local and is applied only to elements that meet explicit admissibility conditions. Similar to refinement, the coarsening procedure is governed by three core principles: (i) hierarchy-aware reversibility, ensuring that only valid parent–child configurations are collapsed; (ii) solution-driven admissibility, based on element-wise indicators derived from field variables; (iii) local conformity preservation, achieved by restricting coarsening to elements whose neighbors satisfy compatibility conditions.

Specifically, coarsening is performed only on sister pairs ($K_1$ and $K_2$) produced by longest-edge bisection from a common mother element $K_m$. A sister pair is admissible for coarsening only if both children are active, *i.e.,* neither child has active descendants, and both satisfy the element-wise coarsening criterion

$$|\eta_K| \geq u_{thr} \tag{48}$$

This criterion identifies regions where is far from the evolving interface and fine resolution is no longer required. To prevent excessive mesh oscillation, coarsening is further restricted to elements above a prescribed minimum level (*i.e.*, $n_{level} > 0$), ensuring that the mesh is never coarsened below its initial resolution.

Local conformity is then enforced by checking the immediate neighbors of $K_1$ and $K_2$: the merge $K_1 \cup K_2 \rightarrow K_m$ is allowed only if it does not create a midpoint (introduced by refinement) that remains on an edge adjacent to a finer element. Practically, this means that all elements that share the bisected edges must be as fine as the post-merge configuration. If any neighbor across the longest edge (or any element incident to the midpoint) is refined to a higher level than the would-be mother element, coarsening is suppressed.

Once a sister pair satisfies all the criteria, they are marked subjective to coarsening. For each admissible sister pair ($K_1$ and $K_2$), coarsening is performed by reactivating the mother element ($K_m$) and deactivating both children. The refinement level of $K_m$ is restored accordingly. The nodal structure is updated by removing mid-edge nodes that are no longer referenced by any active element.

The complete coarsening procedure can be summarized as Algorithm 2 in Figure 7.

**Algorithm 2:** Mesh coarsening *via* hierarchical sister-pair merging

```
Input: Element matrix T, node coordinates r, nodal values u
Output: Updated T, r, u
Compute element-wise indicators η_K for all active elements;
Mark coarsening candidates:
 P ← (|η_K| ≥ u_thr) ∧ is_active(T) ∧ T.level > 0 ∧ is_sister_pair(T) ∧ is_comformity(T) ;
while any(P) do
    /* Admissibility and conformity checks                                  */
    Identify admissible sister pairs (K_1, K_2) with common mother K_m;
    Ensure both children are active leaves at the same refinement level;
    Verify that merging does not violate conformity with neighboring elements;
    /* Merge operation                                                       */
    Deactivate children (K_1, K_2) and reactivate the mother element K_m;
    Update refinement levels, child/sister links, and hierarchy flags;
    Remove unused mid-edge node(s) from r and corresponding values from u;
    Update the candidate set P for the next coarsening pass;
Update local neighbor connectivity;
```

Figure 7 Pseudocode for hierarchical adaptive mesh coarsening based on longest-edge bisection.

*4.5 Dependence of adaptive mesh quality*

This subsection investigates the performance of the proposed HAMR algorithm in generating adaptive meshes that dynamically track the evolving interface in phase-field solution.

Generally, accurately capturing curvature and interface motion in phase-field modeling requires at least three to four grid layers across the diffuse interface. Using the analytical interface thickness presented in Eq. (9), the finest mesh resolution required to resolve the interface properly is given by $\Delta h_{fin} = \frac{5.18\sqrt{\kappa}}{N_I}$, where $N_I$ denotes the desired number of grid layers across the smooth interface. As an example, consider a unit square domain with $\kappa = 0.0004$. If choosing five layers of grid to describe the smooth interface, *i.e.*, $N_I = 5$, we get that the finest mesh resolution would be $\Delta h_{fin} \approx 0.02$, implying that at least a $50 \times 50$ mesh is required. If the coarse mesh is $10 \times 10$, the maximum refinement level, $\bar{n}_{level}$ must be at least 5 ~ 6 to achieve this finest resolution. To validate this estimation, a series of numerical experiments were conducted with the maximum refinement level varying from 2 to 8 while all other parameters were held fixed, as shown in Figure 8. The results confirm that the maximum refinement level controls the number of effective grid layers across the interface: small values (*e.g.*, $\bar{n}_{level} = 2$ in Figure 8a) provide only two layers of refined mesh, which is inadequate for quantitative accuracy. As $\bar{n}_{level}$ increases, more layers of refined mesh appear across the interface, improving accuracy but also increasing the element count. Beyond a certain point, however, additional refinement leads to diminishing accuracy gains while substantially increasing the number of elements and, consequently, the computational effort. Furthermore, to maintain high order accuracy of the GSM solver, the refined layers must remain symmetric; otherwise, accuracy degrades. The numerical results in Figure 8(c) show that $\bar{n}_{level} = 6$ produces a desired mesh with five layers of symmetric and uniform fine mesh across interface, consistent with the analytical estimate and sufficient for the interface thicknesses considered. This alignment confirms that, for any given $\kappa$, the required mesh resolution can be determined automatically based on the relationship in Eq. (9) without manual tuning.

Complementary to the maximum refinement level, the refinement threshold $u_{thr}$ determines how broadly the framework defines the interface zone. The sensitivity of the mesh to $u_{thr}$ is illustrated in Figure 9. Results reveal that a small $u_{thr}$ restricts refinement to only the steepest part of the interface, producing a very narrow band of fine elements. This is computationally efficient but may under-resolve the gradual transition from bulk to interface. As $u_{thr}$ increases, the refined band widens and more elements are included around the interface, improving accuracy while increasing the total element count. However, if $u_{thr}$ becomes too large, the refined zone grows unnecessarily wide and the additional computational cost no longer brings meaningful accuracy gains. In practice, we find that $u_{thr} = 0.925$ provides a good balance between accuracy and efficiency for the phase field applications.

In addition to producing robust static mesh configurations, the adaptive meshing algorithm also supports dynamic remeshing as the interface evolves, as shown in Figure 10, with a circular interface moving from left to right. The refined elements consistently follow the interface, while the bulk remains discretized by the coarsest mesh. The transition zone forms a narrow annular band that moves with the interface, and the mesh hierarchy is updated only where local indicators change. So, remeshing remains localized near the moving interface and the majority of the domain is left untouched. This locality is crucial for computational efficiency, especially in long-time Cahn–Hilliard simulations where interface motion is slow relative to the time-stepping rate. To prevent excessive overhead from frequent remeshing, the update frequency is also regulated. Because the interface moves only a small fraction of a grid cell per time step, updating the mesh every step is unnecessary. Numerical experiments show that updating the mesh every 100 time steps for C-H modeling and 40 time steps for A-C modeling with the suggested values of $\bar{n}_{level}$ and $u_{thr}$ has negligible impact on accuracy while greatly reducing the overall remeshing cost.

Taken together, the hierarchical mesh construction, the automated parameters $\bar{n}_{level}$ and $u_{thr}$, and the controlled remeshing frequency produce an automated, adaptive structured mesh that tightly tracks the evolving phase-field interface, preserves the desired interfacial resolution, and achieves high computational efficiency.

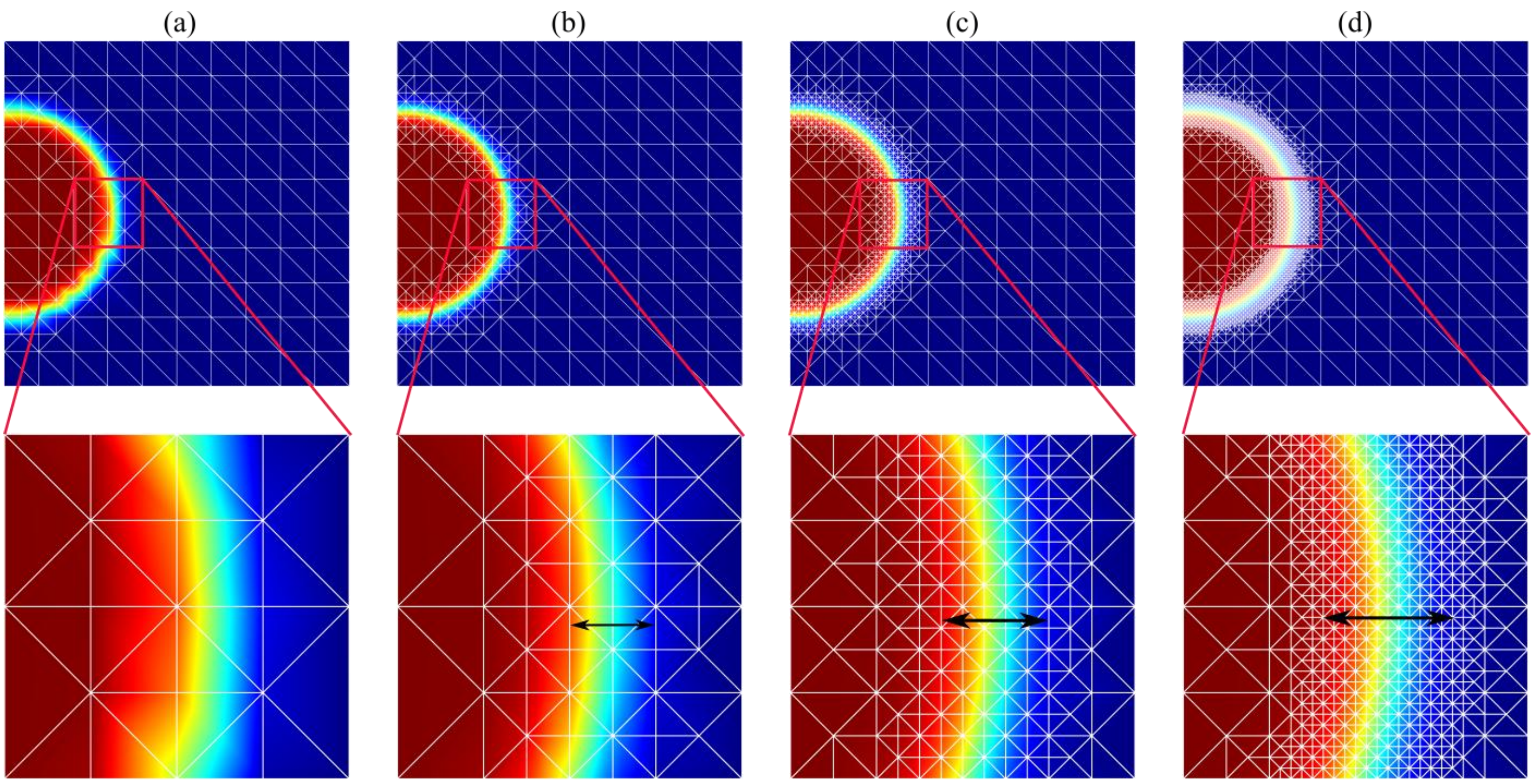

Figure 8 Effect of refinement level n-level on the adaptive mesh resolution near the diffuse interface. Panels (a)–(d) correspond to $\bar{n}_{level}$ = 2, 4, 6, and 8, respectively, with $u_{thr}$ = 0.925 and $\kappa$ = 0.0004. The top row shows the global mesh and concentration field, while the bottom row provides zoom-in views of the interfacial region. The black arrow indicates the region where the nodes have a symmetric neighbors distribution.

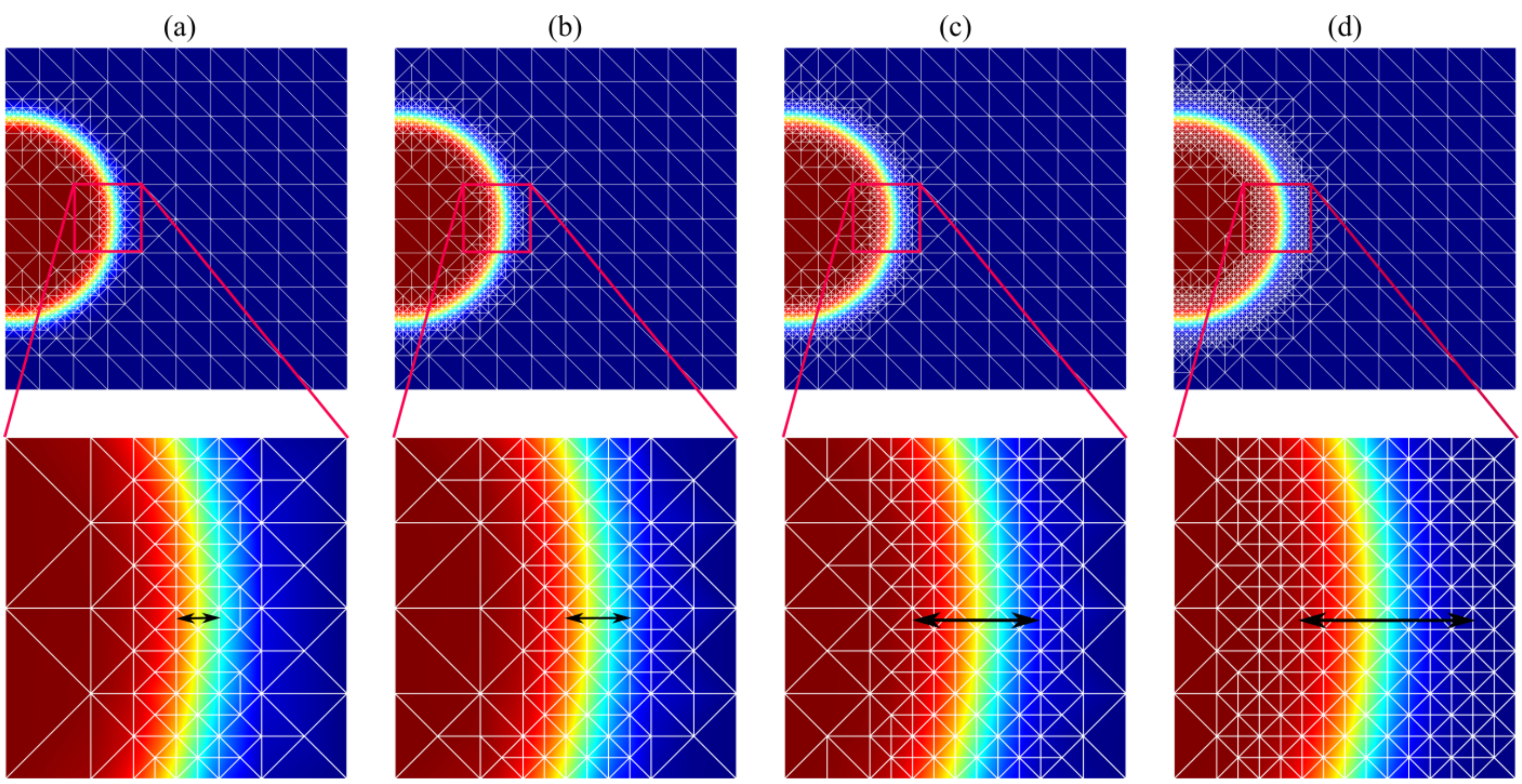

Figure 9 Effect of the refinement threshold $u_thr$ on the quality of the adaptive mesh in the vicinity of the diffuse interface. Panels (a)–(d) correspond to $u_{thr}$ = 0.6, 0.8, 0.925, 0.99, respectively, with $\bar{n}_{level}$ = 6 and $\kappa$ = 0.0004.

The top row shows the global mesh and concentration field, while the bottom row provides zoom-in views of the interfacial region. The black arrow indicates the region where the nodes have a symmetric neighbors distribution.

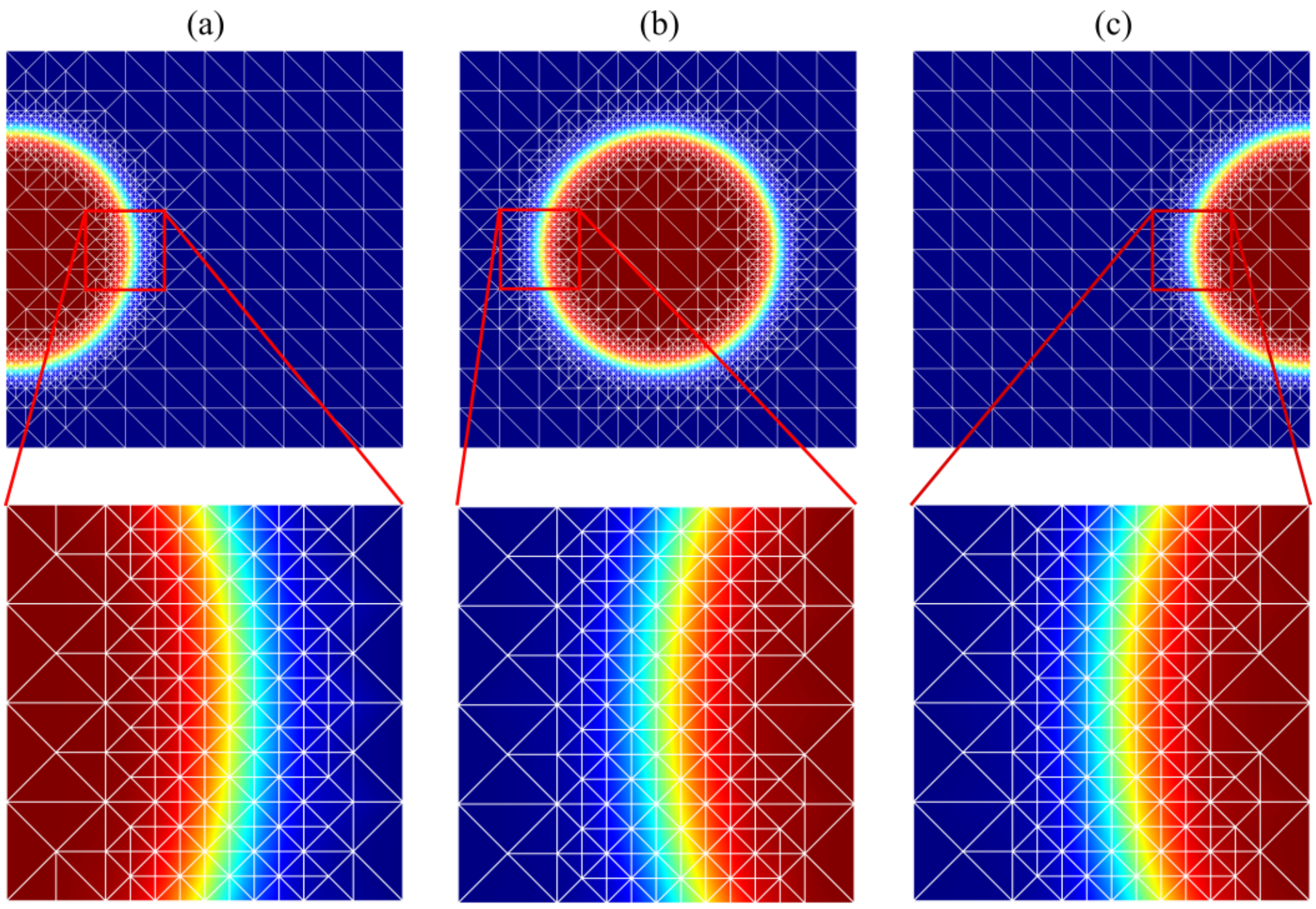

Figure 10 Evolution of the adaptive structured mesh for the GSM solver as the circular interface moves from left to right.

## 5. Results

### *5.1 GSM Laplace solver on uniform mesh*

The performance of the GSM solver for approximating second-derivative terms is first assessed on uniform meshes within the domain of a unit square and compared directly with a standard central-difference finite difference method (FDM), using both accuracy and efficiency as metrics. The test set comprises four smooth functions: a polynomial function

$$F_1(x, y) = x^3 + x^2y^3 \tag{49}$$

a trigonometric function

$$F_2(x, y) = \sin\,(x^2 + y^3) \tag{50}$$

an exponential function

$$F_3(x, y) = e^{x^2+y^3} \tag{51}$$

and a rational function

$$F_4(x, y) = \frac{1}{x^3 + y^2 + 1} \tag{52}$$

For each function, the exact derivative is known analytically, allowing a rigid assessment of the approximation error.

Figure 11 presents the root-mean-squared error (RMSE) of the Laplacian approximation as a function of the grid spacing $h$, where the derivatives are computed by GSM and FDM on successively refined uniform meshes. As shown in Figure 11, the RMSE of GSM exhibits a clear second-order convergence rate with respect to $h$ for all test functions. This confirms that on uniform meshes the GSM achieves the same formal second-order spatial accuracy as the classical central-difference scheme. In particular, the RMSE levels for GSM and FDM are of the same order of magnitude for each resolution, indicating that GSM is quantitatively equivalent to the second-order FDM for Laplacian approximation on uniform meshes. Consequently, the same level of accuracy is expected when GSM is used to discretize the Laplacian terms in phase-field formulations based on uniform mesh.

The computational efficiency of GSM relative to FDM on uniform meshes is summarized in Figure 12, where the CPU time required for a single evaluation of the Laplacian operator is plotted against the total number of grid points $N \times N$. Both implementations are fully vectorized and loop-free in the core operators to ensure a fair comparison. For both methods, the computational cost grows approximately linearly with the number of unknowns, indicating that GSM inherits the same asymptotic computational complexity as FDM when applied on a uniform structured mesh. The results also show that the ratio of CPU times between GSM and FDM remains essentially constant across all resolutions considered: a single GSM operator evaluation is about eight times more expensive than its FDM counterpart. This overhead is primarily attributable to the more involved stencil construction and data access patterns in GSM. In addition, while the FDM Laplacian uses a four-point stencil, the GSM Laplacian involves six neighboring nodes (as shown in Figure 13), which naturally increases the number of floating-point operations per degree of freedom. Consequently, in the purely uniform-mesh setting, FDM is more efficient than GSM in computation and memory storage. The overall performance gain of GSM in efficiency does not stem from superior efficiency on uniform grids. Instead, it results from the reduced number of degrees of freedom afforded by adaptive, non-uniform meshes, which offsets the higher per-element cost and yields substantial savings in both CPU time and memory for large-scale phase-field simulations.

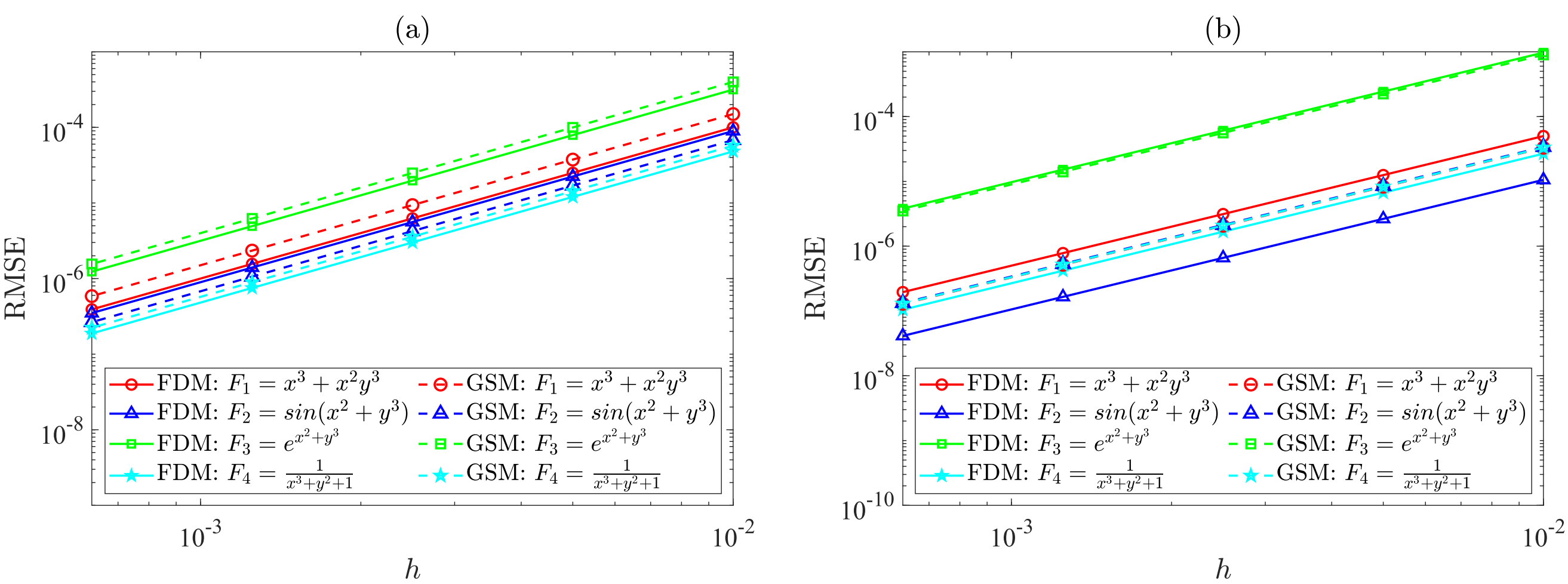

Figure 11 Convergence of the Laplacian approximation on uniform meshes using the GSM and FDM. (a) RMSE of $\partial^2 F/\partial x^2$. (b) RMSE of $\partial^2 F/\partial y^2$.

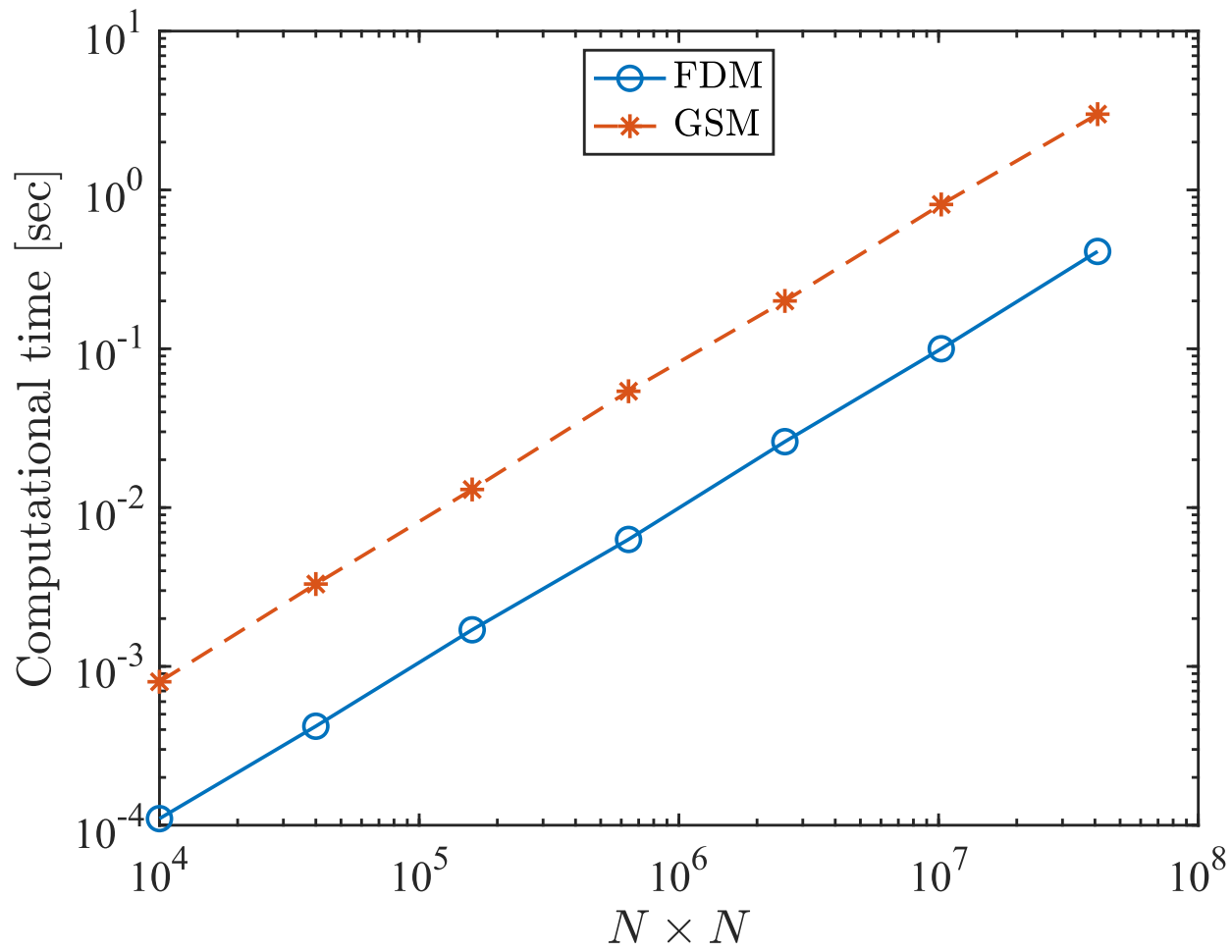

Figure 12 Computational cost of a single Laplacian operator evaluation for GSM and FDM on uniform meshes. The computational time corresponds to the CPU wall-clock time required to apply Laplacian operator evaluation once over the entire domain on a desktop workstation equipped with an intel i9-10940X processor and 64GB RAM.

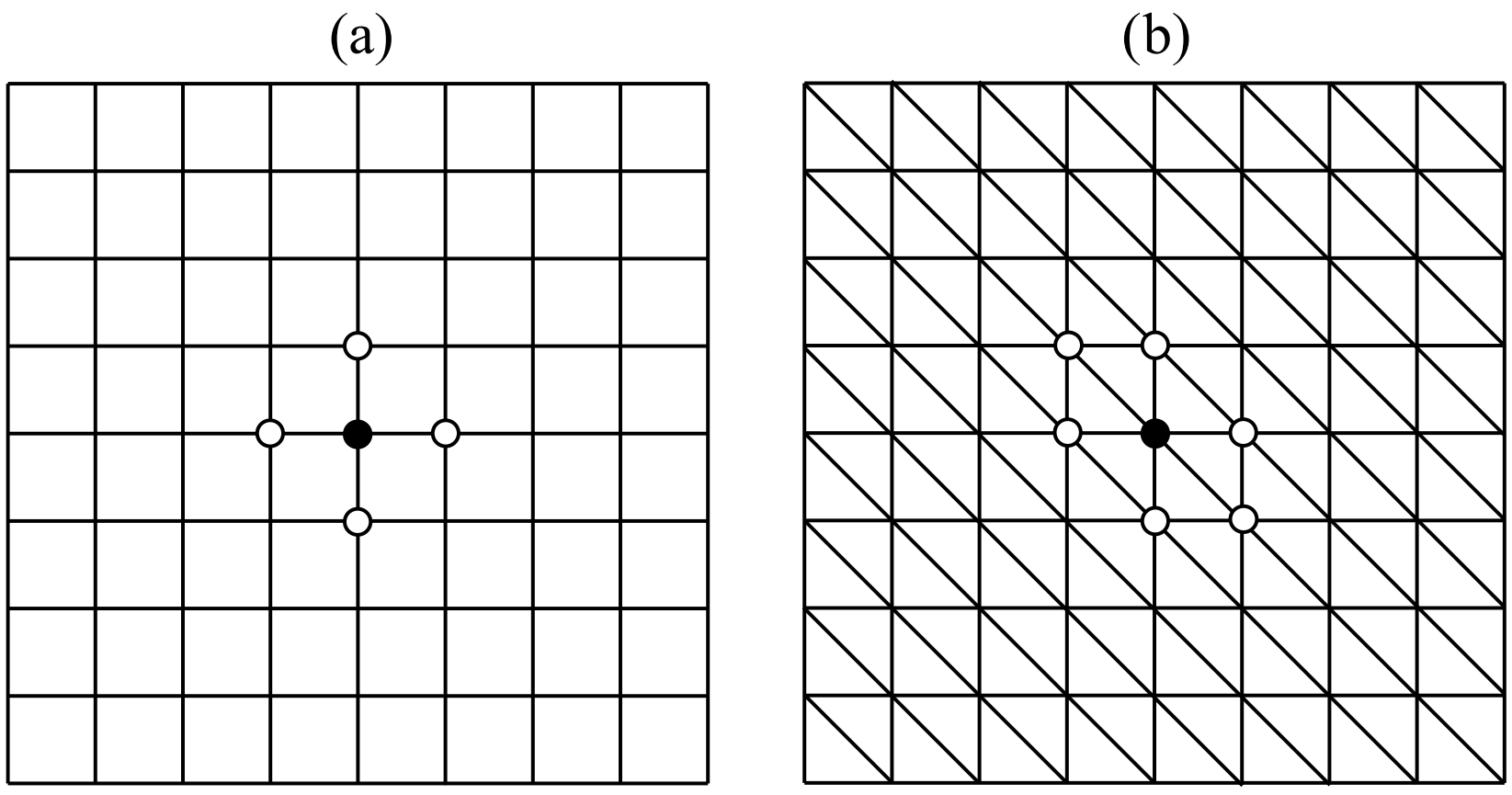

Figure 13 Stencil structure for Laplacian approximation on a uniform Cartesian mesh. (a) Standard FDM stencil using four immediate neighbors on a structured grid. (b) GSM stencil on the corresponding structured triangular mesh, involving six neighboring nodes.

### *5.2 GSM Laplace solver on non-uniform mesh*

The preceding subsection established that, on uniform meshes, the GSM discretization of second-derivative terms is quantitatively equivalent to the standard second-order central-difference FDM in terms of accuracy. In practice, however, the main motivation for using GSM is its ability to accommodate non-uniform and adaptive meshes, in which the stencil structure around each node becomes irregular and generally unsymmetric. In this subsection we examine the accuracy of GSM under such non-uniform conditions and assess how the order of accuracy is affected by local mesh non-uniformity.

To quantify this effect, we consider the Laplacian of the smooth test function $f(x, y) = x^3 + x^2y^3$ evaluated at the point (1, 1). Six representative stencil patterns arising in adaptive mesh, labeled GSM-a through GSM-f in Figure 14, are constructed by introducing extra neighboring nodes while maintaining the same nominal mesh size $h$. Cases (a) and (f) correspond to two uniform stencils with six and eight neighboring nodes, respectively, whereas the remaining cases represent non-uniform stencils with either six or seven neighboring nodes.

Figure 15 presents the convergence of the RMSE as a function of $h$ for these six stencils and compared to the reference FDM on uniform mesh. The results demonstrate that GSM retains second-order accuracy for stencils that remain effectively uniform (GSM-a and GSM-f), consistent with the findings in the preceding subsection. When the stencil becomes appreciably non-uniform (GSM-b–e), the convergence rate degrades to first order, as indicated by the slope of the corresponding error curves in Figure 15. This behavior is consistent with the theoretical expectations for gradient-smoothing schemes on irregular point distributions: the lack of symmetry and uniform spacing introduces additional truncation terms of order $O(h)$ [24, 26].

This observation has direct implications for the design of the adaptive structured mesh employed in the phase-field examples. In the proposed HAMR strategy, both the bulk regions and the diffuse interface are discretized by uniform meshes—a coarse uniform mesh in the bulk and a fine uniform mesh in the interface—so that the GSM discretization in the interface regions preserves second-order accuracy. The non-uniform mesh appears only in a narrow transition band between bulk and interface, where elements of different sizes are connected. From the perspective of the phase-field variables, this transition zone lies within the bulk, where the order parameter or concentration is nearly constant and the physical gradients are very small. As a result, the absolute error introduced by the locally first-order GSM approximation in the non-uniform band is negligible to the overall accuracy of the GSM solver based on adaptive mesh when applied to phase field modeling.

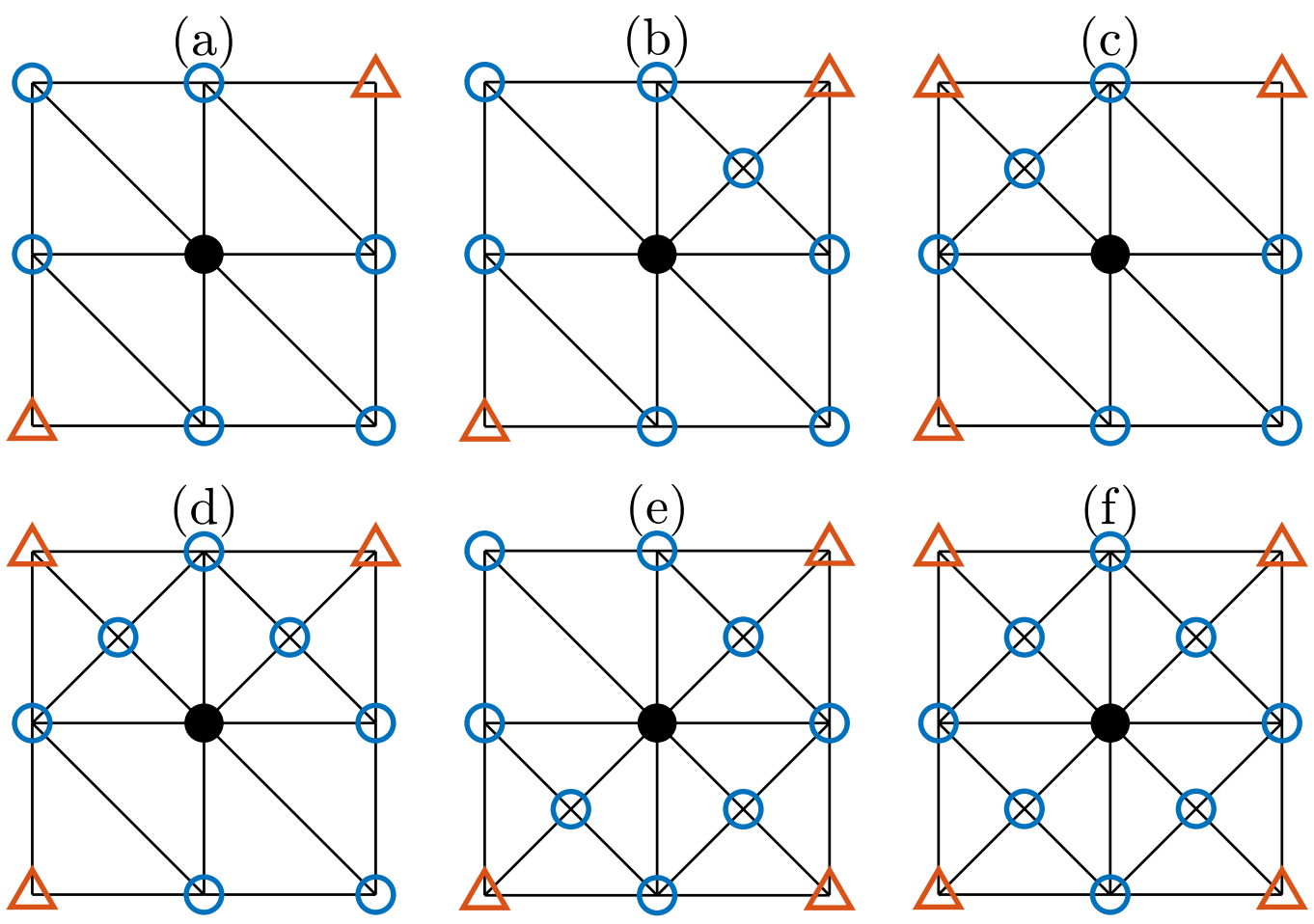

Figure 14 Representative stencil configurations that arise in the adaptive structured mesh for Laplacian approximation using GSM. The black dot denotes the target node, blue circles indicate neighboring nodes that are topologically connected to the target and participate in the stencil, and red triangles mark nearby nodes that are not connected and therefore do not contribute to the derivative approximation.

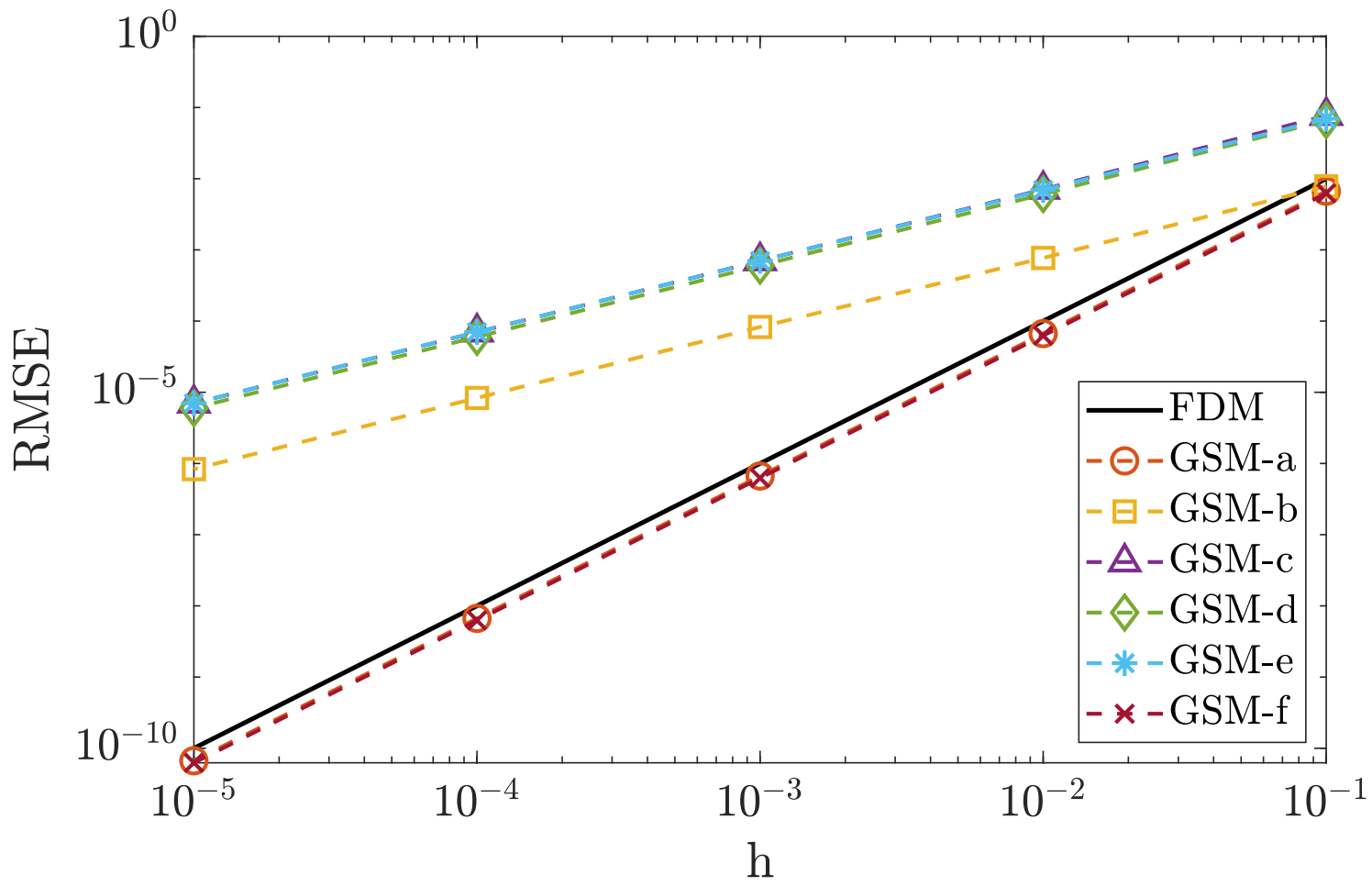

Figure 15 Convergence of the Laplacian approximation on the selected six GSM stencils compared with GDM.

### *5.3 Application to smoothing of an initially sharp interface*

To evaluate the accuracy of the proposed adaptive GSM solver for phase-field modeling, a benchmark problem with an analytical equilibrium solution is first considered. The computational domain is a unit square, in which a sharp vertical interface is initially located at $x = 0.5$, as shown in Figure 16(a). The Allen–Cahn and Cahn–Hilliard equations are then solved using the GSM solver with $\kappa = 0.001$ on both uniform and adaptive meshes. The numerical solutions are compared with the corresponding FDM solutions. At equilibrium, the phase-field variable admits an analytical smooth interface profile, as illustrated in Figure 16(b), defined by the expression in Eq. (8). Figure 16(c) compares the order-parameter distributions along a horizontal line obtained from GSM, FDM, and the analytical solution. Both numerical solvers capture the interface profile accurately.

To quantitatively evaluate the convergence behavior of the numerical solvers with respect to grid size $h$, a series of numerical tests were conducted using different grid resolutions. The errors obtained from the GSM solver on both uniform and adaptive meshes are compared with those from the FDM solver on a uniform mesh, as shown in Figure 17. The results exhibit a clear and consistent convergence trend as the grid is refined. On uniform meshes, the GSM solver demonstrates the same second-order accuracy as the FDM solver, and even slightly better accuracy, which can be attributed to the larger number of supporting neighbors and thus richer local information employed in the derivative approximation. More importantly, the GSM solver based on adaptive meshes also achieves second-order accuracy. Although its error is slightly higher than that of the GSM solver on uniform meshes, it remains smaller than that of the FDM solver. It is worth noting that the GSM solver theoretically attains second-order accuracy only on uniform meshes and degrades to first-order accuracy on non-uniform meshes. However, the present results indicate that the adaptive GSM solver still preserves overall second-order convergence. This improvement is attributed to the adaptive mesh strategy, which maintains completely uniform and fine meshes in the interfacial region where solution evolution is significant, while allowing mesh non-uniformity in the bulk regions close to interface where solution variations are minor and the induced errors are thus negligible. These numerical results confirm that the proposed HAMR algorithm effectively restores the global second-order accuracy of the GSM solver to phase field modeling.

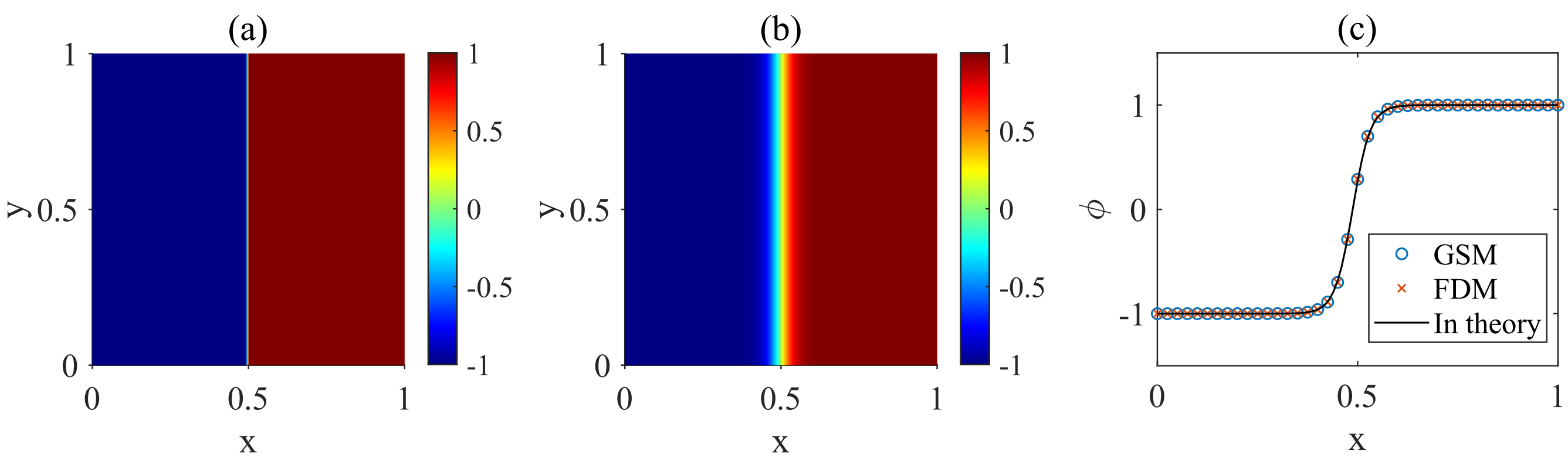

Figure 16 Phase-field benchmark with $\kappa = 0.001$. (a) Initial order-parameter field with a sharp vertical interface at $x = 0.5$. (b) Analytical equilibrium order-parameter field exhibiting a smooth interface. (c) Comparison of order-parameter profiles along a horizontal line obtained from GSM and FDM with the analytical solution in Eq. (8).

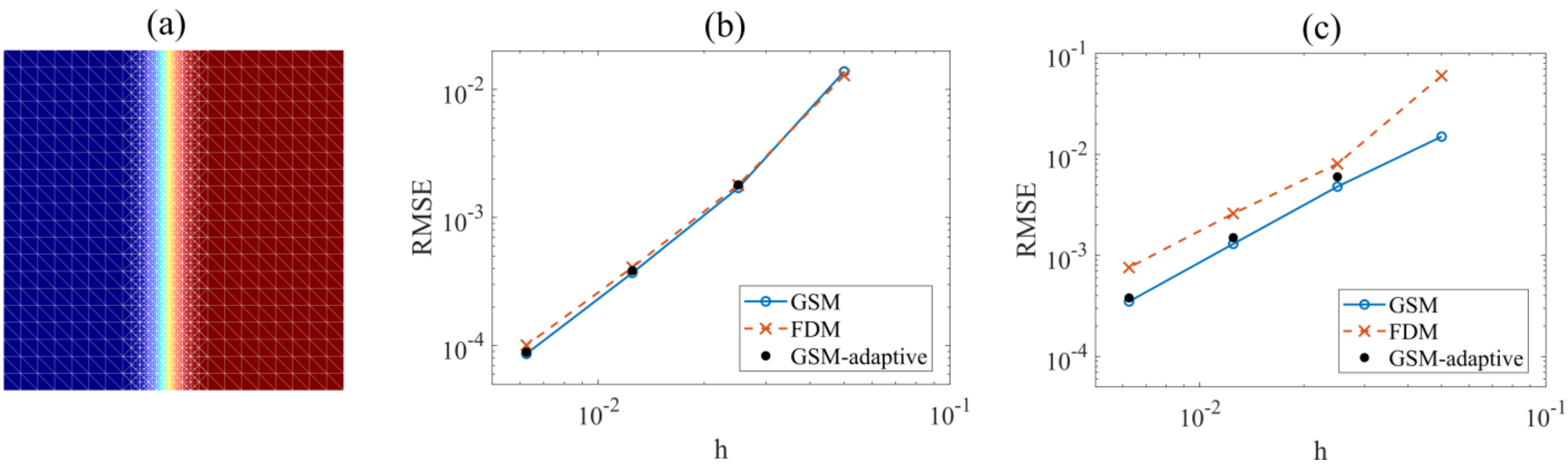

Figure 17 (a) An adaptive mesh sample used in GSM solver with $N_0 = 21$ and $n_{level_{max}} = 4$. (b) Error comparison for the A-C equation with $\kappa = 0.01$ at $t = 1.0$. (c) Error comparison for the C-H equation with $\kappa = 0.001$ at $t = 0.025$. In (b) and (c), GSM and GSM-adaptive correspond to GSM solution on uniform and adaptive meshes, respectively.

*5.4 Application to curvature-driven coarsening*

Next, we assess the accuracy and computational complexity of the proposed adaptive GSM framework for large-scale, long-time phase field modeling by simulating the curvature-driven coarsening. To this end, we consider two representative problems: (i) the evolution of an initially square phase embedded in another phase governed by the A-C equation, and (ii) the evolution of an initially rectangular phase embedded in another phase governed by the C-H equation. In both cases, we vary the interface thickness through the parameter $\kappa$, which directly influences the mesh resolution required to accurately resolve the diffuse interface.

The efficiency of the GSM solver with adaptive mesh refinement are evaluated against the finite difference method (FDM) on a uniform mesh and the FEM-based MOOSE software on adaptive quadrilateral mesh. The five tested conditions (I-V) are summarized in Table 1. Figure 18 illustrates the dependence of interface thickness on $\kappa$ for three representative A-C cases (I, III, and V), along with the corresponding uniform meshes used in FDM and the adaptive meshes used in GSM and MOOSE. In case V, where $\kappa$ is very small, the interface is extremely thin and demands a very fine spatial resolution. For FDM, this necessitates a uniform 1281 × 1281 mesh over the entire domain. In contrast, the GSM framework restricts this fine resolution to a narrow band around the interface, while employing a much coarser mesh away from it. This targeted refinement is the primary source of GSM's advantage in computational efficiency over FDM. An analogous relationship between interface thickness and mesh requirements is also observed for the C-H simulations.

Figure 19 compares interface profiles during the evolution obtained from GSM with adaptive mesh and FDM with uniform mesh. The results obtained from GSM and FDM agree exactly for both the A-C and C-H models, confirming the accuracy of the proposed GSM framework. In the A-C case as shown in Figure 19(a), the embedded square phase shrinks continuously under curvature-driven motion and eventually vanishes, consistent with the non-conserved nature of the order parameter. In the C-H case in Figure 19(b), where the concentrations of the two phases are conserved, the initial rectangular phase evolves toward a circular shape, which minimizes interfacial area and thus interfacial energy, representing the thermodynamically stable state.

Notably, Figure 20 compares the order parameter fields obtained from GSM without and with the directional correction in Eq. (21), highlighting its necessity for accurate Laplacian approximation in the GSM solver. Without this correction, the solution develops a spurious "checkerboard" pattern, as shown in Figure 20(a). When the directional correction is applied, this numerical artifact is completely eliminated, and the interface remains smooth and physically consistent, as demonstrated in Figure 20(b).

Table 1 Five test cases with different values of $\kappa$ and the corresponding mesh parameters. For the adaptive mesh used in GSM, the coarsest mesh in the bulk region is fixed at 21 × 21, while the mesh within the interfacial region, defined by $u \in [-u_{thr}, u_{thr}]$, is refined to the deepest level specified by the $\bar{n}_{level}$.

| **Case #** | $\kappa$ | $h$ | **Uniform mesh** $N \times N$ | **Adaptive mesh** | | |
|---|---|---|---|---|---|---|
| | | | | $N_0$ | $c_{thr}$ | $\bar{n}_{level}$ |
| I | 0.0004 | 1/80 | $81 \times 81$ | 21 | 0.925 | 4 |
| II | 0.0001 | 1/160 | $161 \times 161$ | 21 | 0.925 | 6 |
| III | $2.5 \times 10^{-5}$ | 1/320 | $321 \times 321$ | 21 | 0.925 | 8 |
| IV | $6.25 \times 10^{-6}$ | 1/640 | $641 \times 641$ | 21 | 0.925 | 10 |
| V | $1.5625 \times 10^{-7}$ | 1/1280 | $1281 \times 1281$ | 21 | 0.925 | 10 |

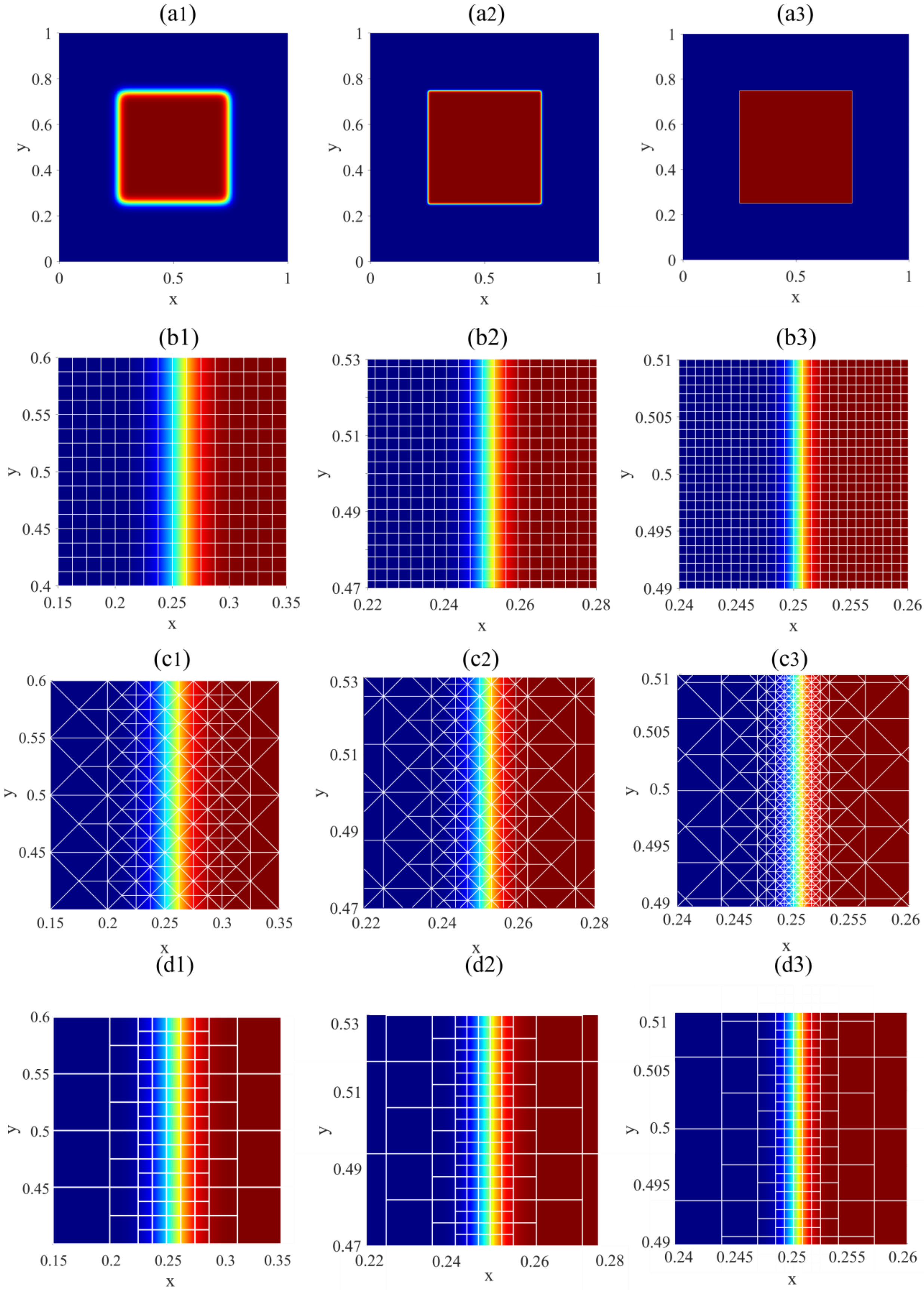

(a1)
(a2)
(a3)
(b1)
(b2)
(b3)
(c1)
(c2)
(c3)
(d1)
(d2)
(d3)
x
y

Figure 18 (a) Initial phase distribution with an equilibrium interface. (b) Uniform mesh used in FDM. (c) Adaptive structured triangular mesh used in GSM. (d) Adaptive structured quadrilateral mesh used in MOOSE. The three columns correspond to the cases I, III, and V in Table 1, respectively.

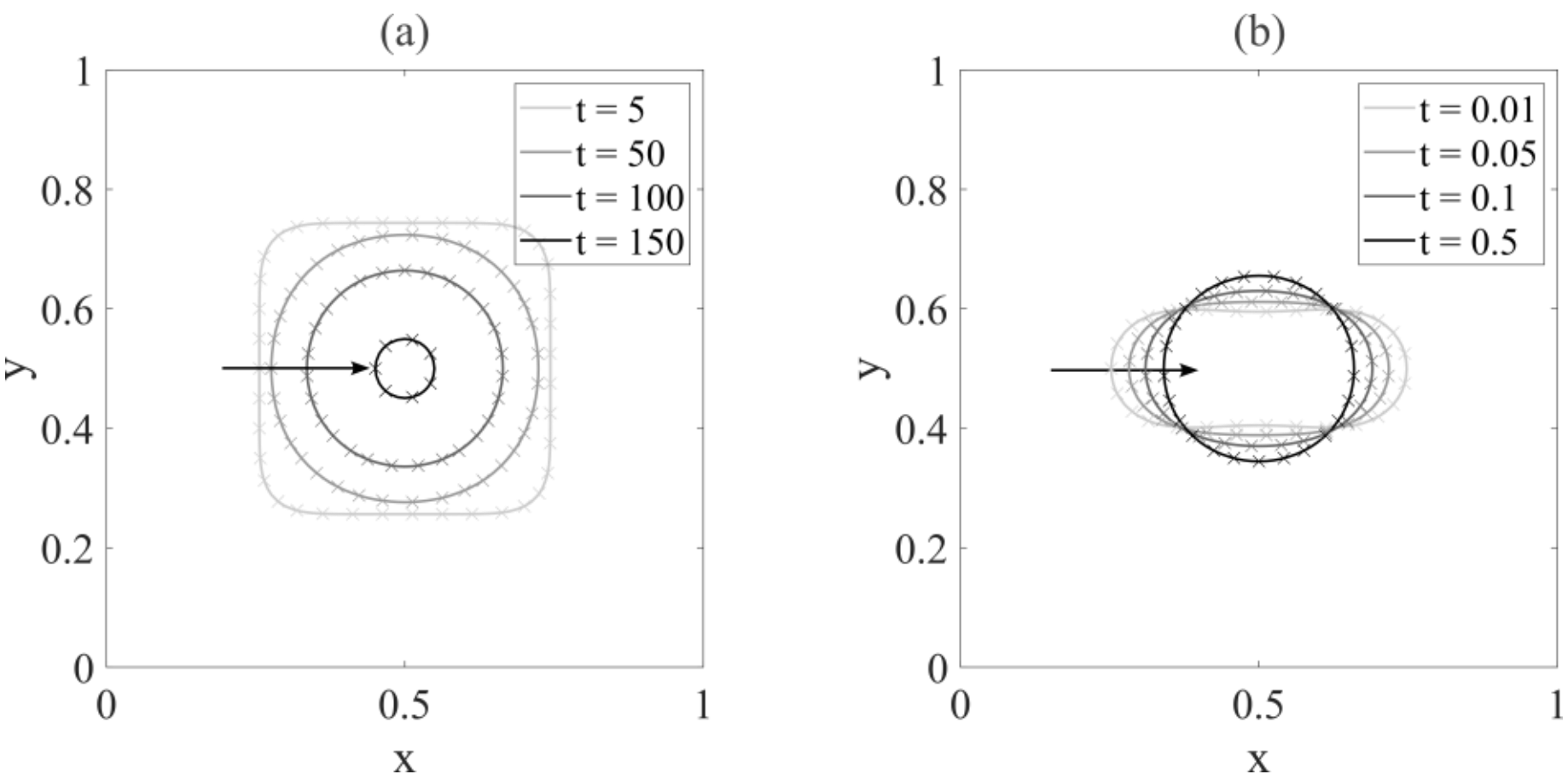

Figure 19 Accuracy comparison between GSM and FDM for (a) A–C and (b) C–H models in case I. The curves show the interface defined by $c = 0$ or $\phi = 0$, evolving from the light-colored shape to the dark-colored shape in the direction indicated by the arrows over time. The solid curves represent FDM results, and the crosses mark the corresponding interfaces obtained from GSM.

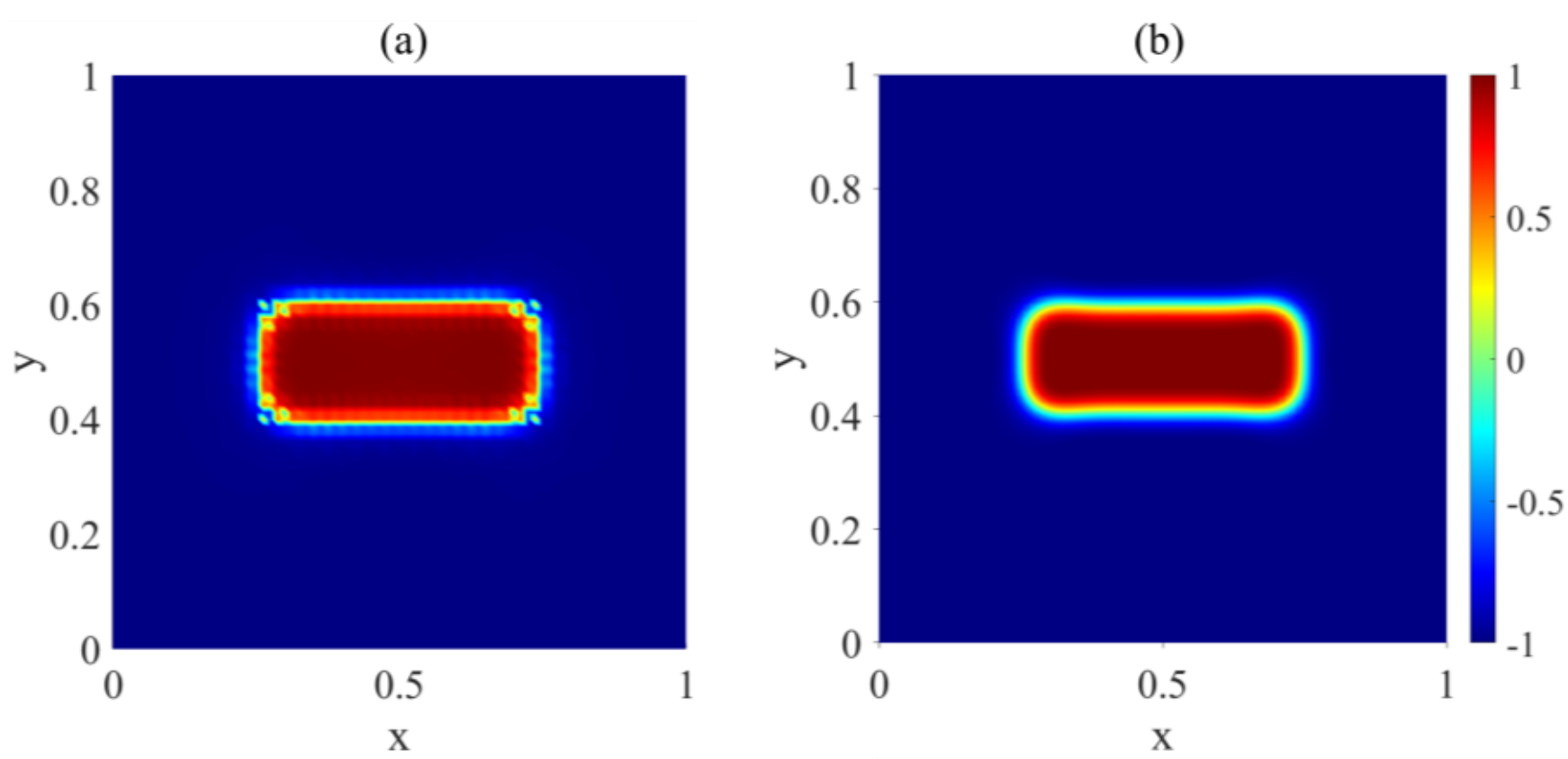

Figure 20 (a) Illustration of the “checkerboard” artifact in the GSM solution without directional correction in C-H modeling of case I. (b) GSM solution with directional correction.

Figure 21 summarizes the computational complexity of the GSM solver with adaptive mesh for C–H modeling and compares it with the standard FDM solver on uniform mesh and FEM solver in MOOSE using adaptive mesh. As shown in Figure 21(a), the FDM runtime increases in proportion to the total number of grid nodes, *i.e.*, $O(N^2)$ computational complexity for a uniform $N \times N$ mesh. In contrast, the GSM solver exhibits a markedly more favorable $O(N)$ scaling because refinement is confined to the diffuse interface region: the computational work grows primarily with the interfacial length rather than with the full domain area. The MOOSE/FEM results also benefit from adaptivity compared with uniform-grid FDM; however, their scaling is less favorable than GSM due to additional overhead from element-wise quadrature, sparse matrix assembly, and iterative solver costs that increase with problem size.

Figure 21(b) further breaks down of computational time among GSM calculation, remeshing, and other tasks in the GSM framework. GSM calculation dominates the total cost, but its fraction decreases with increasing model size $N$, while the fraction of remeshing grows. This trend arises because the GSM calculation step is fully vectorized and thus highly efficient, whereas the remeshing step involves unavoidable loops that cannot be fully vectorized, leading to an overhead that grows comparatively faster in computational cost. In these simulations, remeshing is performed every 100 time-steps, yet its increasing share underlines the importance of further optimizing the meshing algorithm for extremely large-scale problems.

The observed $O(N)$ computational complexity of GSM can be explained by how the number of triangular elements grows with system size $N$. As shown in Figure 22(a), the elements count increases linearly with $N$: along the interface length direction, the number of elements increases with $N$ or decreasing grid size $h$, but in the interface-thickness direction the adaptive mesh always represents the diffuse interface with approximately 4–5 layers of elements, so the element count in this thickness direction remains essentially constant. Consequently, the total number of elements scales as $O(N)$. In contrast, the number of grid points in the FDM mesh grows as $O(N^2)$, consistent with its observed quadratic scaling in CPU time. Figure 22(b) tracks the evolution of the number of triangular elements during a GSM simulation. The total number of elements decreases over time, reflecting the nature of C–H dynamics, which drive the system toward morphologies with reduced interfacial area and lower chemical potential. A smaller interfacial area naturally requires fewer interface-resolving elements. This behavior implies that, beyond the favorable overall $O(N)$ scaling, the effective time-stepping speed of GSM accelerates in the later stages of the simulation as the interface shrinks and the element count decreases.

Comparing GSM with the adaptive FEM solver in MOOSE under the same conditions (grid size, time-step size, and remeshing frequency), GSM remains consistently faster and shows better computational complexity as model size increases. This is because, even though both methods use adaptive meshes, FEM incurs increasing overhead from adaptive mesh management, element-level residual/Jacobian evaluation, sparse matrix assembly, and iterative solver convergence. As a result, the FEM runtime does not scale purely with the element count. Nevertheless, adaptive FEM still becomes more efficient than uniform-grid FDM beyond a sufficiently large problem size, *e.g.,* $700 \times 700$ in this example, reflecting the benefit of concentrating resolution near the interface rather than uniformly across the domain.

Finally, Figure 21(a) indicates that FDM is the most efficient approach for small model sizes. This is expected, since FDM uses a simpler stencil-based Laplacian, whereas GSM requires a more complex Laplacian approximation based on more local particle neighborhoods, which is more time-consuming per degree of freedom. However, due to its lower computational complexity order, the ratio of GSM to FDM runtime decreases as the model size $N$ increases, and GSM becomes faster than FDM once the model size exceeds approximately $400 \times 400$.

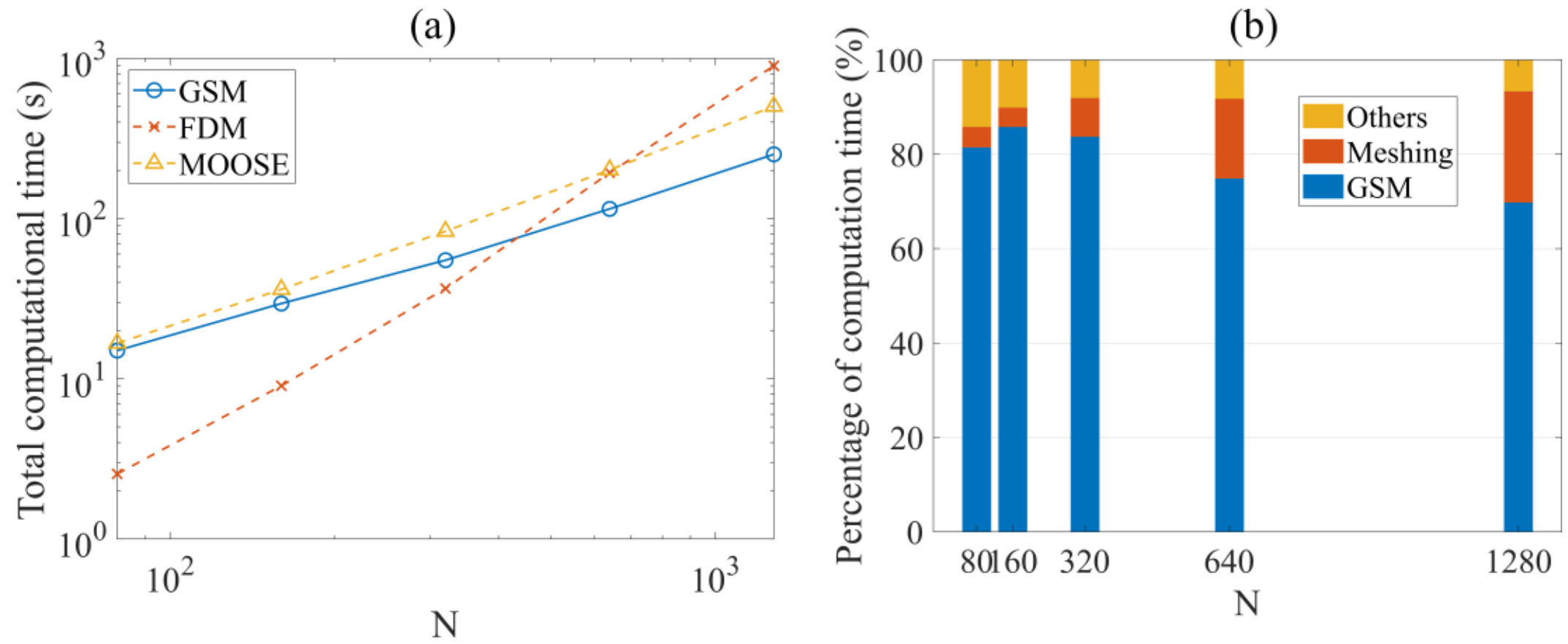

Figure 21 Computational time comparison between GSM, FDM, and FEM-based MOOSE at different spatial resolutions for solving the C–H equation. The reported times correspond to the first 10,000 time steps, obtained on a desktop workstation with an Intel i9-10940X processor and 64 GB RAM. Both GSM and FDM solver are fully vectorized for best efficiency.

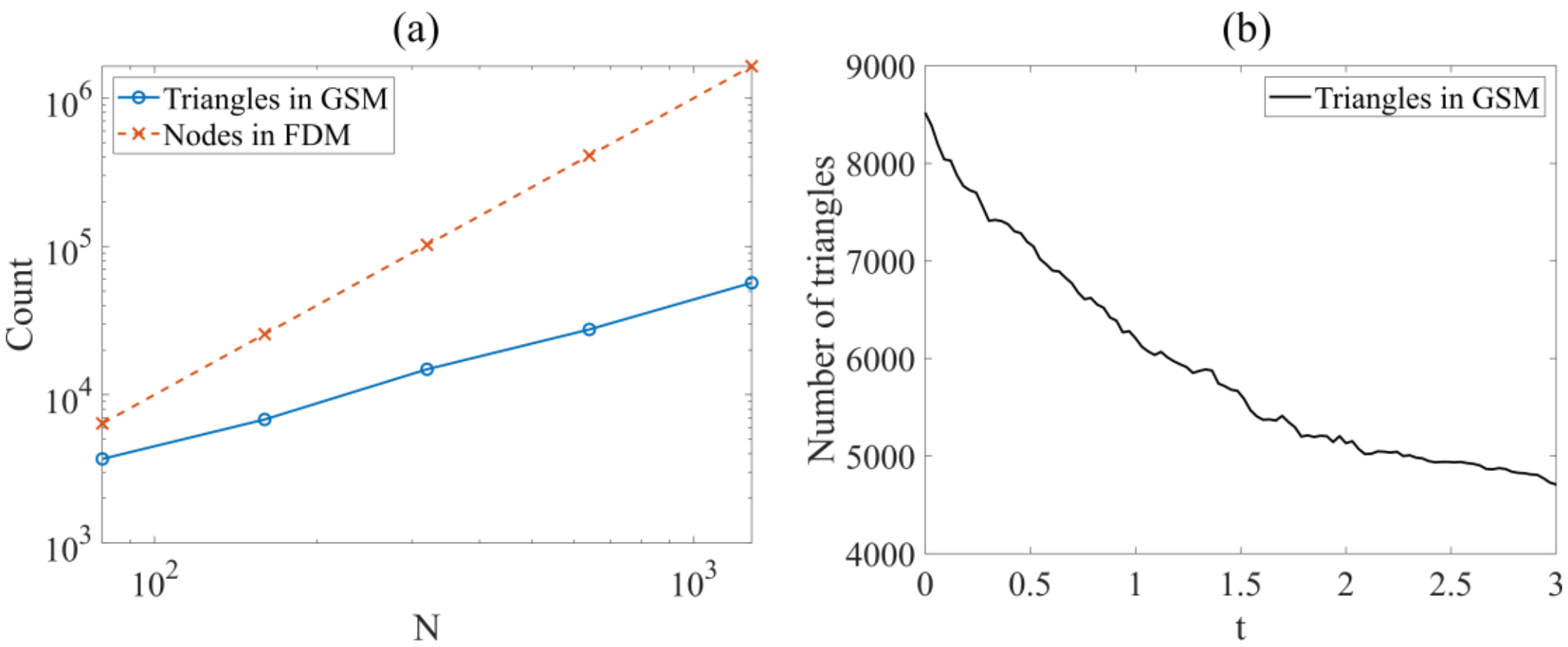

Figure 22 (a) Number of triangular elements in the adaptive mesh for the same interface morphology but different interface thicknesses at varying spatial resolutions. (b) Evolution of the number of triangular elements in the GSM adaptive mesh over the course of the simulation.

The results in Figure 23 for A-C model demonstrate that the numerical solvers for A–C equation exhibits the same efficiency trends as those for C–H equation: the GSM solver scales as $O(N)$, whereas the FDM solver scales as $O(N^2)$. One notable difference is that the fraction of computational time devoted to meshing in the A–C simulations (Figure 23(b)) is higher than in the C–H case (Figure 21(b)) under the same conditions. This can be explained by the different differential orders of the two models. The C–H equation is fourth order in concentration and requires applying the GSM Laplacian operator twice per time step, whereas the A–C equation is second order and uses the Laplacian only once. Consequently, the GSM update step for A–C is intrinsically cheaper, reducing its relative share of the overall runtime and making the meshing cost comparatively more prominent. In addition, the remeshing procedure is performed more frequently in the A–C simulations (every 40 time steps) than in the C–H simulations (every 100 time steps), which further increases the relative contribution of meshing to the overall computational cost in the A–C case.

The distinct physical behaviors of the two models also influence efficiency. In C–H simulations, concentration is conserved, so the interface persists and the system can evolve over long times while maintaining a substantial interfacial area and, therefore, a large number of interface-resolving elements. In A–C simulations, by contrast, the isolated phase A embedded in phase B shrinks monotonically and eventually disappears. Once the embedded phase vanishes, the mesh becomes fully coarse and the number of triangular elements in GSM is maximally reduced, yielding very high efficiency in the late stages. In practice, A–C simulations often terminate after relatively few time steps because the system reaches a stable state quickly, whereas C–H simulations typically require much longer runs and thus impose a heavier computational load. As a result, although GSM still delivers an $O(N)$ scaling advantage for A–C, the practical demand for aggressive acceleration is generally more critical for C–H modeling.

Note that, while the observed $O(N)$ complexity of GSM is broadly applicable to phase-field simulations, the specific crossover point at which GSM becomes more efficient than FDM (around $400 \times 400$ in our tests) is problem dependent. In particular, it depends on the interface morphology, the relative interface area, and the domain size. For other applications with different ratios of interface area to domain size, the system size at which GSM overtakes FDM will change. More rigorously, the crossover is controlled by the ratio between the number of FDM grid nodes and the number of GSM triangular elements. Recall that, as shown in Figure 12, GSM spends roughly eight times more CPU time per triangle than FDM spends per node; thus, once the ratio of FDM nodes to GSM triangles exceeds about eight, GSM becomes the more efficient choice. This condition is easily satisfied in large-scale problems where the interface occupies only a small portion of the domain—either because the interface is short or because it is very thin—since FDM must refine the mesh everywhere, whereas GSM refines only in a narrow band around the interface and keeps the total element count well below $N^2$.

An additional advantage of GSM is reduced memory usage for large models. FDM stores field variables on an $N \times N$ grid, whereas GSM stores them on an adaptive unstructured mesh of size proportional to the number of triangular elements, which is typically much smaller than $N^2$ (see Figure 22(a)). This makes GSM more memory-efficient and thus more suitable for extremely large-scale simulations. Overall, these results indicate that GSM can provide substantial efficiency gains over FDM for very large phase-field models, particularly in cases where the interfacial region occupies only a small fraction of the computational domain, while maintaining at least comparable, and often superior, accuracy.

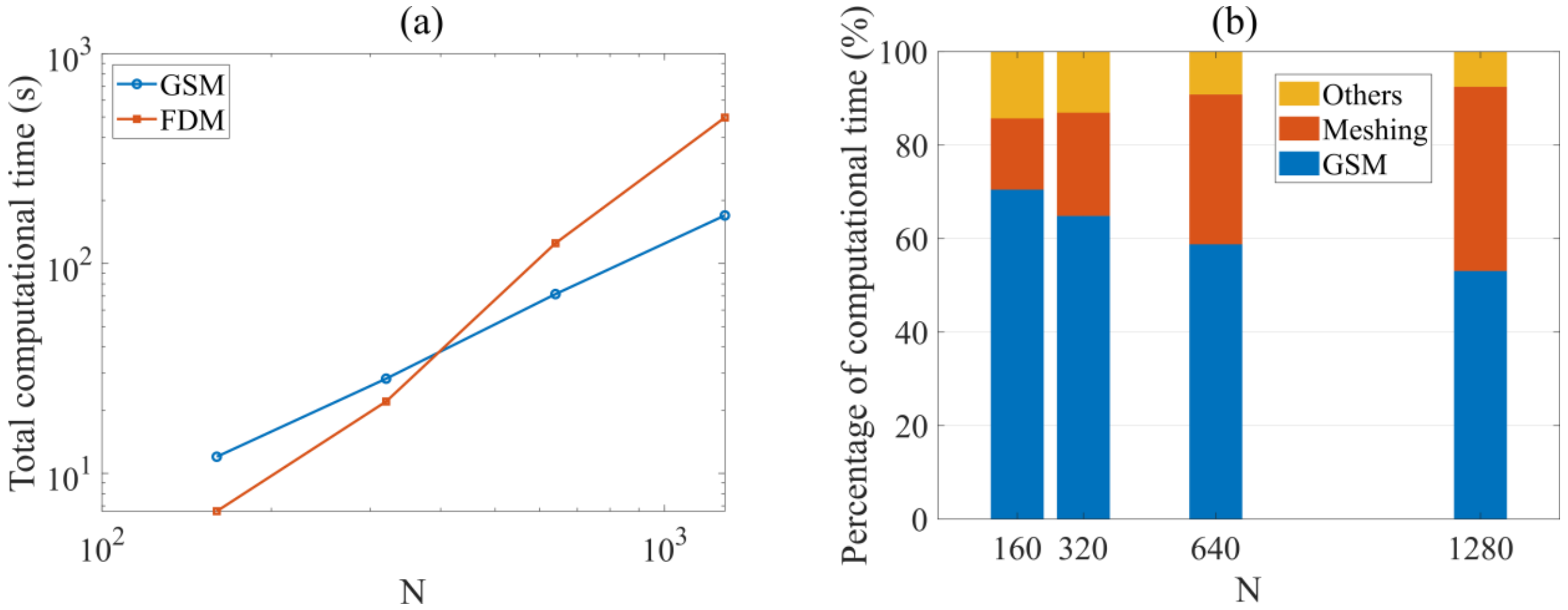

Figure 23 Efficiency comparison between GSM and FDM for solving the A–C equation. The reported times correspond to the first 10,000 time steps, obtained on a desktop workstation with an Intel i9-10940X processor and 64 GB RAM. Both GSM and FDM solver are fully vectorized for best efficiency.

### *5.5 Application to cases with non-uniform mobility*

The preceding phase-field examples assumed a constant, spatially uniform mobility. In many applications, however, the mobility depends on the order parameter (and may vary in space and time). Such non-uniform mobility increases the stiffness and nonlinearity of the governing equations.

To demonstrate that adaptive GSM can effectively handle phase-field models with strongly nonlinear, non-uniform mobility, we simulate the curvature-driven coarsening problem using the Cahn–Hilliard (C–H) model with a concentration-dependent mobility. Specifically, we adopt the mobility $M(c) = M_0|1 - c^2|$ proposed by Cahn *et al.* [43] to represent interface-dominated transport, *i.e.*, diffusion localized primarily within the diffuse interface. We refer to this choice as *interfacial mobility*. Using this mobility, we simulate the annealing and retraction of a terminating layer with the adaptive GSM framework and benchmark the results against the corresponding FDM solver.

Figure 24 compares the retraction profiles predicted by adaptive GSM and FDM under interfacial mobility. The two solutions are visually indistinguishable throughout the evolution, demonstrating the accuracy and robustness of adaptive GSM for phase-field models with strongly nonlinear, non-uniform mobility. The retraction process proceeds as follows: at early times, the initially sharp corner at the terminating edge rapidly rounds due to curvature and the layer begins to retract leftwards. As retraction continues, a bulge develops in the terminating end. At later times, a neck forms ahead of the bulge and progressively thins until pinch-off occurs, resulting in detachment of the bulge from the remaining layer. Notably, this behavior differs qualitatively from the uniform-mobility case, in which the layer tends to retract into a single large droplet rather than undergoing repeated pinch-off. These results indicate that spatially varying mobility can substantially alter the retraction pathway and promote more complex morphological evolution.

We further investigated the influence of interface thickness, controlled by the gradient-energy coefficient $\kappa$ on the retracting dynamics. A series of simulations with different $\kappa$ values is summarized in Figure 25. As $\kappa$ decreases (thinner diffuse interface), the layer breaks into fewer but larger droplets; conversely, increasing $\kappa$ (thicker interface) produces more numerous but smaller droplets. This trend can be understood from the role of $\kappa$ in the phase-field free energy defined in Eq. (1): the gradient term penalizes spatial variations of the order parameter and thus provides an intrinsic regularization that influences which morphological length scales are dynamically accessible during retraction. Since breakup proceeds *via* a capillarity-driven necking instability, the characteristic spacing between necks (and hence between pinch-off events) depends on the interfacial thickness. For larger $\kappa$, the interface is thicker and the order parameter varies more gradually across the interfacial zone. Two consequences follow. First, curvature and the associated chemical-potential gradients are distributed over a broader region, so the capillary driving is less concentrated at a single sharp location and can more readily activate multiple nearby sites along the retracting edge. Second, the wider interfacial zone provides a larger effective region for interfacial transport, enabling small perturbations along the rim to grow and persist rather than rapidly merging into a single dominant bulge. The instability therefore develops at a shorter characteristic spacing, leading to more frequent neck formation and pinch-off and, consequently, a larger number of smaller droplets (top rows of Figure 25). In contrast, for smaller $\kappa$ the interface becomes thinner and the evolution is more localized: fewer necks develop, pinch-off occurs less often, and the same material is partitioned into fewer detached regions. By mass conservation, each droplet then contains more material and attains a larger radius, consistent with the reduced droplet count and increased droplet size observed in the lower rows of Figure 25.

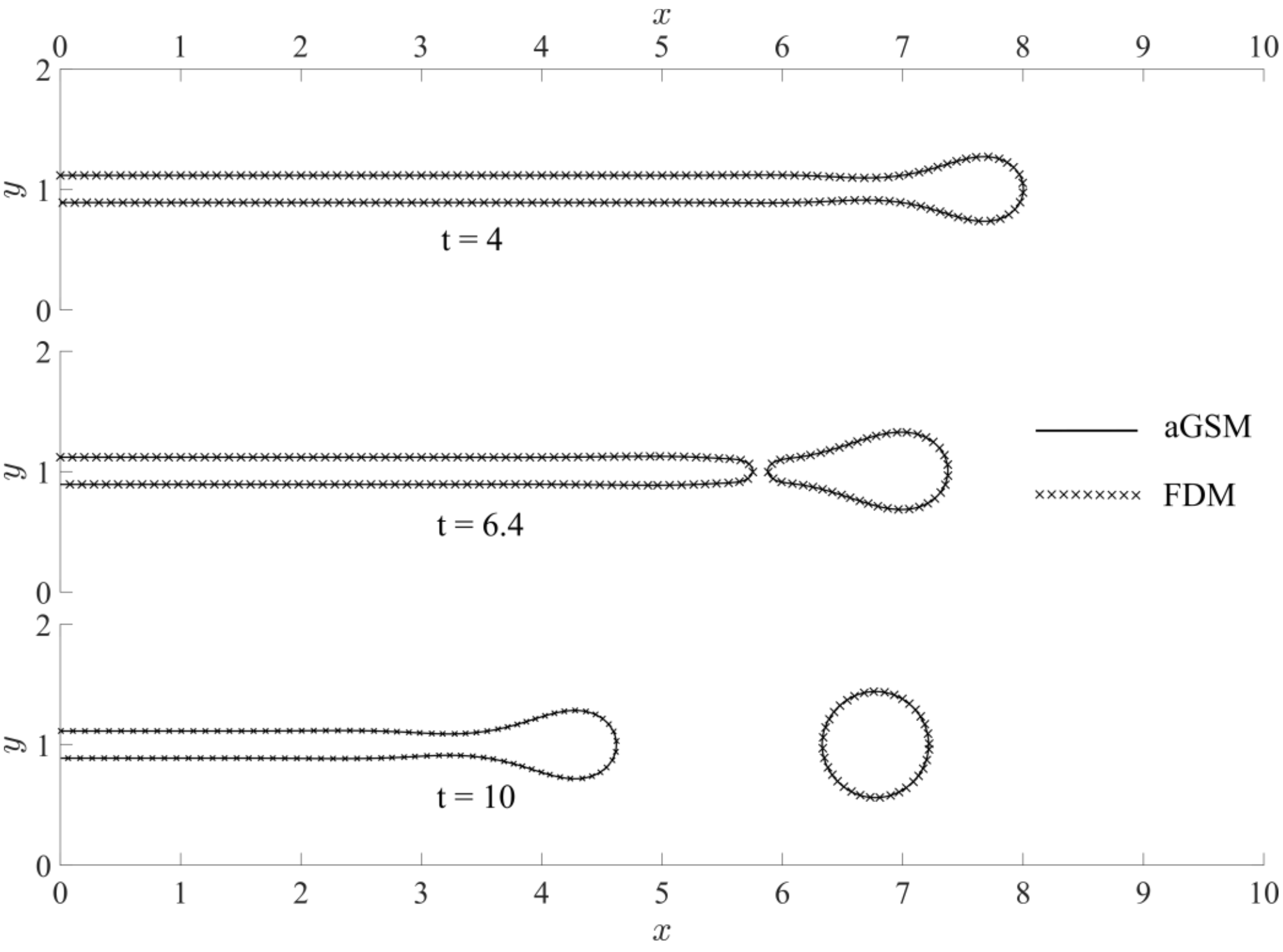

Figure 24 Comparison of interface profiles during layer retraction with interfacial mobility obtained from adaptive GSM (denoted as aGSM) and benchmarking FDM. $M_0 = 1$, $\kappa = 0.0004$. The initial layer size is $9 \times 0.2$.

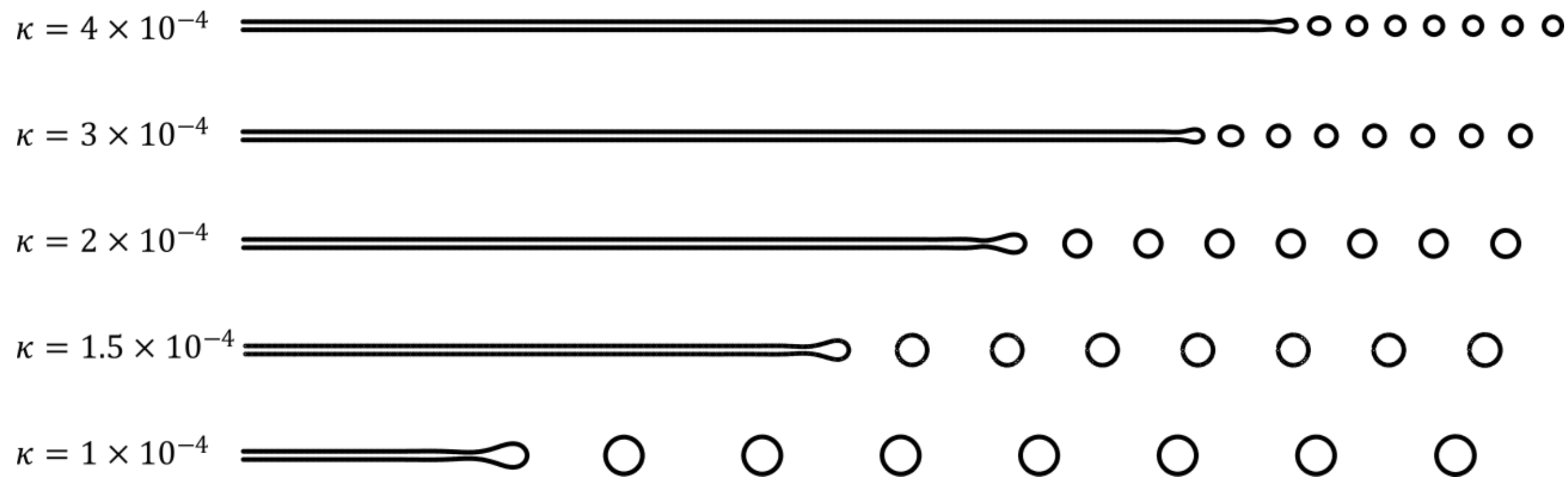

Figure 25 Comparison of interface profiles under different $\kappa$, obtained from adaptive GSM with interfacial mobility. $M_0 = 1$.

## 6. Conclusions and discussion

This work demonstrates that the gradient smoothing method (GSM), combined with a hierarchical adaptive structured triangular mesh, is a practical and efficient alternative to uniform-grid finite difference methods (FDM) for large scale phase-field simulations with thin interface.

In accuracy, the GSM discretization achieves second-order convergence on uniform meshes, matching the second-order accurate FDM and in most cases delivering slightly smaller errors, thanks to the richer local information provided by a larger supporting neighbor set in the derivative approximation. On non-uniform meshes, GSM is theoretically expected to degrade to first-order accuracy; however, the proposed adaptive strategy restores an overall second-order convergence rate in practice for both the Allen–Cahn and Cahn–Hilliard equations. This is because the adaptive algorithm preserves a completely uniform, fine mesh in the interfacial region where gradients and higher-order derivatives dominate the solution evolution, while relegating mesh non-uniformity to bulk regions near the interface where the dynamics is weak and the resulting discretization errors are negligible. The adaptive remeshing procedure can automatically generate high-quality meshes that consistently provide the required number of resolution layers across the diffuse interface and tracking its evolution, ensuring accurate interface representation for both second-order (A–C) and fourth-order (C–H) phase-field dynamics.

In terms of computational efficiency, the GSM solver incurs a higher cost per degree of freedom than the uniform-grid FDM solver, owing to its broader stencil, gradient smoothing operations, and more complex data management. Consequently, for small problems or fully uniform resolutions, FDM can be faster. However, the primary advantage of GSM lies in its compatibility with adaptive meshing. By confining fine resolution to the narrow diffuse interface while maintaining a coarse mesh in the bulk regions, the GSM solver significantly reduces the total number of degrees of freedom. This reduction leads to substantially lower memory requirements and mitigates the per-node overhead of GSM, making it more suitable for large-scale problems involving large domains and/or very thin interfaces. From the computational complexity perspective, GSM exhibits favorable $O(N)$ complexity thanks to its applicability to adaptive meshes, in contrast to the $O(N^2)$ scaling of uniform-grid FDM. That is, when the problem size (indicated by $N$) increases, the overall wall-clock time and memory footprint of GSM grow linearly with $N$ rather than with a globally refined $N \times N$ grid. This characteristic enables GSM to outperform FDM in large-scale phase-field simulations, say $N > 400$ in the tested example.

As a strong-form discretization, the finite-difference method (FDM) is typically much more computationally efficient than general-purpose FEM frameworks such as MOOSE when using the same number of degrees of freedom (DOF) and the same time-step size. This advantage arises because FDM applies Laplacian operators through fixed local stencils with regular memory access (high cache locality) and minimal abstraction overhead, often in a naturally matrix-free manner. In contrast, FEM implementations generally incur higher per-degree-of-freedom costs due to element-wise quadrature, residual and Jacobian assembly, sparse matrix insertion, and the expense of iterative linear/nonlinear solver convergence. Consequently, for comparable DOF counts and time steps, FDM commonly achieves typically one order of magnitude speed advantage per time step. However, classical FDM is not readily extensible to non-uniform, unstructured meshes, limiting its ability to further reduce cost *via* adaptive refinement like FEM solvers. The proposed adaptive GSM framework retains the efficiency advantages of strong-form, stencil-like operators while enabling adaptive meshing, thereby achieving efficiency gains similar to adaptive FEM, but with substantially lower operator and assembly overhead than weak-form FEM solvers. Therefore, under the same conditions, the adaptive GSM framework is expected to be much more efficient than both FDM and FEM solvers for large-scale problems.

Overall, the results demonstrate that, while retaining FDM-level second-order accuracy, the GSM framework coupled with adaptive meshing provides superior computational complexity and efficiency for large-scale phase-field applications.

Several promising directions can be pursued to further enhance the efficiency, computational complexity, and applicability of the GSM framework. While both FDM and FEM benefit from mature implicit or semi-implicit time-integration schemes that permit much larger stable time steps, the present GSM implementation uses explicit time stepping, which imposes restrictive stability limits for stiff phase-field models, particularly the Cahn–Hilliard equation. Extending the adaptive GSM framework to implicit or semi-implicit schemes is thus necessary. Notably, the stability criteria and theoretical foundations applied to the strong-form FDM solver are directly applicable to the strong-form GSM framework as well. However, efficient numerical algorithms must be specifically developed to make these schemes fully compatible with the GSM discretization and data structures.

The GSM formulation and the hierarchical adaptive mesh refinement strategy are inherently dimension-independent and can be naturally extended to three-dimensional applications. In three dimensions, GSM based on adaptive mesh is expected to exhibit $O(N^2)$ complexity, which remains significantly more favorable than the $O(N^3)$ complexity of uniform-grid FDM. However, extending to large-scale 3D simulations entails substantial computational cost and long runtimes. Fortunately, both the GSM derivative operator and the adaptive meshing procedure operate on a node-wise basis using only local information. This intrinsic locality makes the GSM framework highly amenable to parallel implementation. In future work, GPU acceleration and large-scale parallel computing will be investigated to further enhance computational performance. By mapping GSM operators and mesh adaptation procedures onto massively parallel architectures, significant reductions in wall-clock time and promising scalability can be expected for large-scale simulations with hundreds of millions of adaptive elements.

**Declaration of Competing Interest**

The authors declare that they have no known competing financial interests or personal relationships that could have appeared to influence the work reported in this paper.

**CRediT authorship contribution statement**

**Zirui Mao**: Investigation, Conceptualization, Methodology, Visualization, Funding acquisition, Resource, Writing - original draft. **Alice Xu**: Methodology, Programming code transfer, Discussion. **Shenyang Hu**: Discussion, Writing - review & editing. **Panos Stinis**: Discussion, Writing – review & editing. **Yuyan Shao**: Discussion, Writing – review & editing. **Yulan Li**: Discussion, Writing – review & editing. **Zhijie (Jay) Xu**: Discussion, Writing – review & editing.

**Acknowledgements**

The work described in this article was performed by Pacific Northwest National Laboratory (PNNL), which is operated by Battelle for the U.S. Department of Energy under Contract DE-AC05-76RL01830. This research was fully supported by the Laboratory Directed Research and Development (LDRD) program at PNNL.